%% file: otsenkimodulei.tex
\documentclass[11pt,reqno]{amsart}
\usepackage{amssymb,amscd,amsbsy,mathrsfs}

\input{classstartscr}

\theoremstyle{remark}

\newtheorem*{rem*}{Remark}

\newcommand{\Bbbone}{{\rm{1\mathchoice{\kern-0.25em}{\kern-0.25em}{\kern-0.2em}{\kern-0.2em}I}}}

\newcommand\Li{{\rm Lip}}
\newcommand\fM{\frak M}

\newcommand\dg{\frak D}

\newcommand\Abs{\operatorname{Abs}}
\newcommand\mX{\mathcal{X}}
\newcommand\mY{\mathcal{Y}}
\newcommand\mI{\mathcal{I}}
\newcommand\mB{\mathcal{B}}
\newcommand\fF{\frak F}
\newcommand\card{\operatorname{card}}
\newcommand\fd{\frak d}
\newcommand\mP{\mathcal{P}}

\begin{document}

\newcommand{\vse}{\vspace{.2in}}
\numberwithin{equation}{section}

\title{Estimates of operator moduli of continuity}
\author{A.B. Aleksandrov and V.V. Peller}
\thanks{The first author is partially supported by RFBR grant 11-01-00526-a ; the second author is partially supported by NSF grant DMS 1001844}

\newcommand{\mt}{{\mathcal T}}

\begin{abstract}
In \cite{AP2} we obtained general estimates of the operator moduli of continuity of functions on the real line. In this paper we improve the estimates obtained in \cite{AP2} for certain special classes of functions.

In particular, we improve estimates of Kato \cite{Ka} and show that
$$
\big\|\,|S|-|T|\,\big\|\le C\|S-T\|\log\left(2+\log\frac{\|S\|+\|T\|}{\|S-T\|}\right)
$$
for every bounded operators $S$ and $T$ on Hilbert space.
Here $|S|\df(S^*S)^{1/2}$. Moreover, we show that this inequality is sharp.

We prove in this paper that if $f$ is a nondecreasing continuous function on $\R$ that vanishes on $(-\be,0]$ and is concave on $[0,\be)$, then its operator modulus of continuity
$\O_f$ admits the estimate
$$
\O_f(\d)\le\const\int_e^\be\frac{f(\d t)\,dt}{t^2\log t},\quad\d>0.
$$

We also study the problem of sharpness of estimates obtained in \cite{AP2} and \cite{AP4}. We construct a $C^\be$ function $f$ on $\R$ such that $\|f\|_{L^\be}\le1$, $\|f\|_{\Li}\le1$, and
$$
\O_f(\d)\ge\const\,\d\sqrt{\log\frac2\d},\quad\d\in(0,1].
$$

In the last section of the paper we obtain sharp estimates of $\|f(A)-f(B)\|$ in the case when the spectrum of $A$ has $n$ points. Moreover, we obtain a more general result in terms of the $\e$-entropy of the spectrum  that also improves the estimate of the operator moduli of continuity of Lipschitz functions on finite intervals, which was obtained in \cite{AP2}.

\end{abstract}

\maketitle

\

\begin{center}
{\Large Contents}
\end{center}

\

\begin{enumerate}
\item[1.] Introduction \quad\dotfill \pageref{intr}
\item[2.]  Schur multipliers \quad\dotfill \pageref{mult}
\item[3.] Remarks on absolutely convergent Fourier integrals  \quad\dotfill \pageref{FouIn}
\item[4.] Estimates of certain multiplier norms \quad\dotfill \pageref{Estim}
\item[5.] Operator Lipschitz functions and operator modulus of continuity \quad\dotfill \pageref{OL}
\item[6.] The operator Lipschitz norm of the function $|x|$
on subsets of $\R$ \quad\dotfill \pageref{|x|}
\item[7.] The operator modulus of continuity of a certain piecewise linear function \quad\dotfill \pageref{kuslin}
\item[8.] Operator moduli of continuity of concave functions on $\R_+$ \quad\dotfill \pageref{conc}
\item[9.] Lower estimates for operator moduli of continuity \quad\dotfill \pageref{lest}
\item[10.] Lower estimates in the case of unitary operators \quad\dotfill \pageref{lest+}
\item[11.] Self-adjoint operators with finite spectrum. Estimates in terms of the $\e$-entropy of the spectrum \quad\dotfill \pageref{safinsp}
\item[] References \quad\dotfill \pageref{bibl}
\end{enumerate}

\

\section{\bf Introduction}
\setcounter{equation}{0}
\label{intr}

\

In this paper we study operator moduli of continuity of functions on subsets of the real line. For a closed subset $\fF$ of the real line $\R$ and for a continuous function $f$ on $\fF$, the {\it operator modulus of continuity} $\O_{f,\fF}$ is defined by
$$
\O_{f,\fF}(\d)\df\sup\|f(A)-f(B)\|,\quad\d>0,
$$
where the supremum is taken over all self-adjoint operators $A$ and $B$ such that
$$
\s(A)\subset\fF,\quad\s(B)\subset\fF,\quad\mbox{and}\quad\|A-B\|\le\d.
$$
If $\fF=\R$, we use the notation $\O_f\df\O_{f,\R}$.
Recall that a continuous function $f$ on $\fF$ is called {\it operator Lipschitz} if
$\O_{f,\fF}(\d)\le\const\d$, $\d>0$.

It turns out that a {\it Lipschitz function} $f$ on $\R$, i.e., a function $f$ satisfying
$$
|f(x)-f(y)|\le\const|x-y|,\quad x,~y\in\R,
$$
does not have to be operator Lipschitz. This was established for the first time by Farforovskaya \cite{F1}. It was shown later in \cite{Ka} that the function
$x\mapsto|x|$ on $\R$ is not operator Lipschitz. The paper \cite{Ka} followed the paper \cite{Mc}, in which it was shown that the function $x\mapsto|x|$ is not commutator Lipschitz. We refer the reader to \S\,\ref{OL} for the definition of commutator Lipschitz functions. Note that nowadays it is well known that operator Lipschitzness is equivalent to commutator Lipschitzness.

We would like to also mention that in \cite{Pe1} necessary conditions for operator Lipschitzness were found that also imply that Lipschitzness is not sufficient for operator Lipschitzness. On the other hand, it was shown in \cite{Pe1} and \cite{Pe3} that if $f$ belongs to the Besov class $B_{\be1}^1(\R)$, then $f$ is operator Lipschitz (we refer the reader to \cite{Pee} and \cite{Pe4} for the definition of Besov classes).

In our joint papers \cite{AP1} and \cite{AP2} we obtain the following upper estimate for continuous functions $f$ on $\R$:
\bay
\label{modnep}
\O_f(\d)\le\const\,\d\int_\d^\be\frac{\o_f(t)}{t^2}\,dt=\const\int_1^\infty\frac{\o_f(t\d)}{t^{2}}ds,\quad\d>0,
\ey
where $\o_f$ is the modulus of continuity of $f$, i.e.,
$$
\o_f(\d)\df\sup\big\{|f(x)-f(y)|:~x,\,y\in\R,~|x-y|\le\d\big\},\quad\d>0.
$$
We deduced from \rf{modnep} in \cite{AP2} that for a Lipschitz function $f$ on $[a,b]$,
the following estimate for the operator modulus of continuity $\O_f$ holds:
$$
\O_{f,[a,b]}(\d)\le\const\,\d\left(1+\log\left(\frac{b-a}{\d}\right)\right)\|f\|_{\Li},
$$
where
$$
\|f\|_{\Li}\df\sup_{x\ne y}\frac{|f(x)-f(y)|}{|x-y|}.
$$
A similar estimate was obtained earlier in \cite{Ka} in the very special case $f(x)=|x|$. Namely, it was shown in \cite{Ka} that for bounded self-adjoint operators $A$ and $B$ on Hilbert space, the following inequality holds:
$$
\big\|\,|A|-|B|\,\big\|\le\frac2\pi\|A-B\|\left(2+\log\frac{\|A\|+\|B\|}{\|A-B\|}\right).
$$
It turns out, however, that for the function $x\mapsto|x|$ the operator modulus of continuity admits a much better estimate. Namely, we show in \S\,\ref{|x|} that
under the same hypotheses
$$
\big\|\,|A|-|B|\,\big\|\le\const\|A-B\|\log\left(2+\log\frac{\|A\|+\|B\|}{\|A-B\|}\right).
$$
We also prove in this paper that this estimate is sharp.

Note that in \cite{NF} an estimate slightly weaker than \rf{modnep} was obtained by a different method.

In \S\,\ref{conc} we show that if $f$ is a continuous nondecreasing function on $\R$ such that $f(x)=0$ for $x\le0$ and the restriction of $f$ to $[0,\be)$ is a concave function, then estimate \rf{modnep} can also be improved considerably:
$$
\O_f(\d)\le\const\int_e^\be\frac{f(\d t)\,dt}{t^2\log t},\quad\d>0.
$$
We also obtain other estimates of operator moduli of continuity in \S\,\ref{conc}.

It is still unknown whether inequality \rf{modnep} is sharp.
It follows easily from \rf{modnep} that if $f$ is a function on $\R$ such that $\|f\|_{L^\be}\le1$, $\|f\|_{\Li}\le1$, then
$$
\O_f(\d)\le\const\,\d\left(1+\log\frac1\d\right),\quad\d\in(0,1].
$$
We construct in \S\,\ref{lest}  a $C^\be$ function $f$ on $\R$ such that $\|f\|_{L^\be}\le1$, $\|f\|_{\Li}\le1$, and
$$
\O_f(\d)\ge\const\,\d\sqrt{\log\frac2\d},\quad\d\in(0,1].
$$
To construct such a function $f$, we use necessary conditions for operator Lipschitzness found in \cite{Pe1}.
We do not know whether the results of \S\,\ref{lest} are sharp.

In \S\,\ref{lest+} we obtain lower estimates in the case of functions on the unit circle and unitary operators.

Finally, we obtain in \S\ref{safinsp} the following sharp estimate for the norms $\|f(A)-f(B)\|$ for Lipschitz functions $f$ and self-adjoint operators $A$ and $B$ on Hilbert space such that the spectrum $\s(A)$ of $A$ has $n$ points:
\bay
\label{ntochek}
\|f(A)-f(B)\|\le C(1+\log n)\|f\|_{\rm Lip}\|A-B\|.
\ey
Moreover, we obtain in \S\ref{safinsp} an upper estimate in the general case (see Theorem \ref{qcom}) in terms of the $\e$-entropy of the spectrum of $A$, where $\e=\|A-B\|$.
It includes inequalities \rf{modnep} and \rf{ntochek} as special cases. Note that \rf{ntochek} improves earlier estimates in \cite{F1} and \cite{F2}.

In \S\,\ref{mult} we give a brief introduction to Schur multipliers, in \S\,\ref{FouIn} we collect auxiliary estimates of certain functions in the space of functions with absolutely converging Fourier integrals. The estimates collected in \S\,\ref{FouIn} are used in \S\,\ref{Estim} to estimate the Schur multiplier norms of certain functions of two variables.
To obtain upper estimates for operator moduli of continuity of concave functions, we estimate in \S\,\ref{kuslin} the operator modulus of continuity of a very special piecewise continuous function on $\R$.

\

\section{\bf Schur multipliers}
\setcounter{equation}{0}
\label{mult}

\

In this section we define Schur multipliers and discuss their properties. Note that the notion of a Schur multiplier can be defined in the case of two spectral measures (see e.g.,
\cite{Pe1}). In this section we give define Schur multipliers in the case of two scalar measures. This corresponds to the case of spectral measures of multiplicity 1.

Let $(\mathcal X,\mu)$ and $(\mathcal Y,\nu)$ be $\s$-finite measure spaces.
Let $k\in L^2(\mathcal X\times\mathcal Y,\mu\otimes\nu)$. Then $k$ induces
the integral operator $\mathcal I_k=\mathcal I_k^{\mu,\nu}$ from
$L^2(\mathcal Y,\nu)$ to $L^2(\mathcal X,\mu)$ defined by
$$
(\mathcal I_k f)(x)=\int_{\mathcal Y}k(x,y)f(y)\,d\nu(y),\quad f\in L^2(\mathcal Y,\nu).
$$
Denote by $\|k\|_{\mathcal B}=\|k\|_{\mB_{\mX,\mY}^{\,\mu,\nu}}$ the operator norm
of $\mathcal I_k$. Let $\Phi$ be a $\mu\otimes\nu$-measurable function defined almost everywhere on
$\mX\times\mY$.
We say that $\Phi$ is a {\it Schur multiplier with respect to $\mu$ and $\nu$} if
$$
\|\Phi\|_{\frak M_{\mX,\mY}^{\,\mu,\nu}}
\df\sup\big\{\|\Phi k\|_{\mB}:
~k\,,\Phi k\in L^2(\mathcal X\times\mathcal Y,\mu\otimes\nu),~\|k\|_{\mB}\le1\big\}<\be.
$$
We denote by $\frak M_{\mX,\mY}^{\,\mu,\nu}$ the space of Schur multipliers with respect to $\mu$ and $\nu$.
It can be shown easily that $\frak M_{\mX,\mY}^{\,\mu,\nu}\subset L^\be(\mX\times\mY,\mu\otimes\nu)$ and $\|\Phi\|_{L^\be(\mX\times\mY,\mu\otimes\nu)}\le
\|\Phi\|_{\frak M_{\mX,\mY}^{\,\mu,\nu}}$.
Thus if $\Phi\in\frak M_{\mX,\mY}^{\,\mu,\nu}$, then
$$
\|\Phi\|_{\frak M_{\mX,\mY}^{\,\mu,\nu}}
=\sup\big\{\|\Phi k\|_{\mB}:
~k\in L^2(\mathcal X\times\mathcal Y,\mu\otimes\nu),~\|k\|_{\mB}\le1\big\}.
$$
Note that $\frak M_{\mX,\mY}^{\,\mu,\nu}$ is a Banach algebra:
$$
\|\Phi_1\Phi_2\|_{\frak M_{\mX,\mY}^{\,\mu,\nu}}
\le\|\Phi_1\|_{\frak M_{\mX,\mY}^{\,\mu,\nu}}\|\Phi_2\|_{\frak M_{\mX,\mY}^{\,\mu,\nu}}.
$$
It is easy to see that $\|\Phi\|_{\frak M_{\mX,\mY}^{\,\mu,\nu}}=\|\Psi\|_{\frak M_{\mY,\mX}^{\,\nu,\mu}}$
for $\Psi(y,x)=\Phi(x,y)$.

If $\mathcal X_0$ is a $\mu$-measurable subset of $X$, then we denote by $(\mathcal X_0,\mu)$
the corresponding measure space on the $\s$-algebra of $\mu$-measurable subsets of $\mathcal X_0$.

Let $\mathcal X=\bigcup\limits_{n=1}^\be\mathcal X_n$ and $\mathcal Y=\bigcup\limits_{n=1}^\be\mathcal Y_n$,
where the $\mathcal X_n$ are $\mu$-measurable subsets
of $\mathcal X$, and the $\mathcal Y_n$ are $\nu$-measurable subsets
of $\mathcal Y$. It is easy to see that
$$
\sup_{m,n\ge1}\|k\|_{\mB_{\mX_n,\mY_n}^{\,\mu,\nu}}^2
\le\|k\|_{\mB_{\mX,\mY}^{\,\mu,\nu}}^2
\le\sum_{m=1}^\be\sum_{n=1}^\be\|k\|_{\mB_{\mX_n,\mY_n}^{\,\mu,\nu}}^2
$$
for every $k\in L^2(\mathcal X\times\mathcal Y,\mu\otimes\nu)$, and
\bay
\label{ml2}
\sup_{m,n\ge1}\|\Phi\|_{\frak M_{\mX_n,\mY_n}^{\,\mu,\nu}}^2
\le\|\Phi\|_{\frak M_{\mX,\mY}^{\,\mu,\nu}}^2
\le\sum_{m=1}^\be\sum_{n=1}^\be\|\Phi\|_{\frak M_{\mX_n,\mY_n}^{\,\mu,\nu}}^2
\ey
for every $\Phi\in L^\be(\mathcal X\times\mathcal Y,\mu\otimes\nu)$.

We state the following elementary theorem:

\begin{thm}
\label{2mult}
Let $(\mathcal X,\mu)$, $(\mathcal X,\mu_0)$, $(\mathcal Y,\nu)$ and $(\mathcal Y,\nu_0)$ be $\s$-finite measure spaces.
Suppose that $\mu_0$ is absolutely continuous with respect to $\mu$ and $\nu_0$ is absolutely continuous
with respect to $\nu$. Let $\Phi\in\frak M_{\mX,\mY}^{\,\mu,\nu}$.
Then $\Phi\in\frak M_{\mX,\mY}^{\,\mu_0,\nu_0}$ and
$\|\Phi\|_{\frak M_{\mX,\mY}^{\,\mu_0,\nu_0}}\le\|\Phi\|_{\frak M_{\mX,\mY}^{\,\mu,\nu}}$.
\end{thm}
\Pf By the Radon--Nikodym theorem, $d\mu_0=\f d\mu$ and $d\nu_0=\psi d\nu$
for nonnegative measurable functions $\f$ and $\psi$ on $\mX$
and $\mY$. Let $k\in L^2(\mX\times\mY,\mu_0\otimes\nu_0)$. Put
$$
(Tk)(x,y)\df k(x,y)\sqrt{\f(x)\psi(y)}.
$$
Clearly, $T$ is an isometric embedding
from $L^2(\mX\times\mY,\mu_0\otimes\nu_0)$ in $L^2(\mX\times\mY,\mu\otimes\nu)$.
Moreover, $\|Tk\|_{\mB_{\mX,\mY}^{\,\mu,\nu}}=\|k\|_{\mB_{\mX,\mY}^{\,\mu_0,\nu_0}}$.
We have
\begin{align*}
\|\Phi k\|_{\mB_{\mX,\mY}^{\,\mu_0,\nu_0}}=&\|T(\Phi k)\|_{\mB_{\mX,\mY}^{\,\mu,\nu}}
=\|\Phi Tk\|_{\mB_{\mX,\mY}^{\,\mu,\nu}}\\[.2cm]
\le&\|\Phi\|_{\frak M_{\mX,\mY}^{\,\mu,\nu}}
\|Tk\|_{\mB_{\mX,\mY}^{\,\mu,\nu}}
=\|\Phi\|_{\frak M_{\mX,\mY}^{\,\mu,\nu}}
\|k\|_{\mB_{\mX,\mY}^{\,\mu_0,\nu_0}}
\end{align*}
for every $k\in L^2(\mX\times\mY,\mu_0\otimes\nu_0)$. Hence, $\Phi\in\frak M_{\mX,\mY}^{\,\mu_0,\nu_0}$ and
$\|\Phi\|_{\frak M_{\mX,\mY}^{\,\mu_0,\nu_0}}\le\|\Phi\|_{\frak M_{\mX,\mY}^{\,\mu,\nu}}$. $\bl$

Note that if $\mX$ and $\mY$ coincide with the set $\Z_+$ of nonnegative integers and $\mu$ and $\nu$ are the counting measure, the above definition coincides with the definition of Schur multipliers on the space of matrices: a matrix $A=\{a_{jk}\}_{j,k\ge0}$ is called a Schur multiplier on the space of bounded matrices if
$$
A\star B\quad\mbox{is a matrix of a bounded operator, whenever}\quad B\quad\mbox{is}.
$$
Here we use the notation
\bay
\label{ScHad}
A\star B=\{a_{jk}b_{jk}\}_{j,k\ge0}
\ey
for the Schur--Hadamard product of the matrices $A=\{a_{jk}\}_{j,k\ge0}$ and $B=\{b_{jk}\}_{j,k\ge0}$.

Let $\mX$ and $\mY$ be closed subsets of $\R$.
We denote by ${\frak M}_{\mX,\mY}$ the space of
Borel Schur multipliers on $\mX\times\mY$, i.e., the space of
Borel functions $\Phi$ defined everywhere on $\mX\times\mY$ such that
$$
\|\Phi\|_{\frak M_{\mX,\mY}}\df\sup\|\Phi\|_{\frak M_{\mX,\mY}^{\mu,\nu}}<\be,
$$
where the supremum is taken over all regular positive Borel measures $\mu$ and $\nu$ on $\mathcal X$ and $\mathcal Y$. It can be shown easily that
$$
\sup_{(x,y)\in\mX\times\mY}|\Phi(x,y)|\le\|\Phi\|_{\frak M_{\mX,\mY}}.
$$
It is also easy to verify that if $\Phi_n\in\fM_{\mX,\mY}$,
$\Phi$ is a bounded Borel function on $\mX\times\mY$, and
$\Phi_n(x,y)\to \Phi(x,y)$ for all $(x,y)\in\mathcal X\times\mathcal Y$,
then
$$
\|\Phi\|_{\frak M_{\mX,\mY}}
\le\liminf\limits_{n\to\be}
\|\Phi_n\|_{\frak M_{\mX,\mY}}.
$$
In particular, $\Phi\in\frak M_{\mX,\mY}$ if
$\liminf\limits_{n\to\be}
\|\Phi_n\|_{\frak M_{\mX,\mY}}<\be$.

We are going to deal with functions $f$ on $\mX\times\mY$ that are continuous in each variable. It must be a well-known fact that such a function $f$ has to be a Borel function.
Indeed, one can construct an increasing sequence $\{\mY_n\}_{n=1}^\be$
of discrete closed subsets of $\mY$ such that $\bigcup\limits_{n=1}^\be\mY_n$
is dense in $\mY$. Let us consider the function $f_n:\mX\times\R\to\C$
such that $f\big|(\mX\times\mY_n)=f_n\big|(\mX\times\mY_n)$ and $f_n(x,\cdot)$ is
a piecewise linear function with nodes in $\mY_n$ for all $x\in\mX$.
Clearly, the function $f_n$ is defined uniquely if we require
that $f_n(x,\cdot)$ is constant on each unbounded
complimentary interval of $\mY_n$.
It is easy to see that $f_n$ is continuous on $\mX\times\R$ and
$\lim\limits_{n\to\be}f_n(x,y)=f(x,y)$ for all $(x,y)\in\mX\times\mY$.
Thus, $f$ belongs to the first Baire class, and so it is Borel.

\begin{thm}
\label{2mult+}
Let $\mX$ and $\mY$ be closed subsets of $\R$ and let $\Phi$ be a
function on $\mX\times\mY$ that is continuous in each variables.
Suppose that $\mu$ and $\mu_0$ are positive regular Borel measures on $\mX$, and
$\nu$ and $\nu_0$ are positive regular Borel measures on $\mY$.
If $\supp\mu_0\subset\supp\mu$ and $\supp\nu_0\subset\supp\nu$,
then $\|\Phi\|_{\frak M_{\mX,\mY}^{\,\mu_0,\nu_0}}\le\|\Phi\|_{\frak M_{\mX,\mY}^{\,\mu,\nu}}$.
\end{thm}

We need two lemmata.

\begin{lem}
Let $\mX$ and $\mY$ be compact subsets of $\R$ and
let $\mu$ and $\nu$ be finite positive Borel measures on $\mX$ and $\mY$.
Suppose that $\{\nu_j\}_{j=1}^\be$ is a sequence  of finite positive Borel measures on $\mY$ that
converges to $\nu$ in the weak-$*$ topology $\s\big((C(\mY))^*,C(\mY)\big)$.
If $k$ is a bounded Borel function on $\mX\times\mY$ such that
$k(x,\cdot)\in C(\mY)$ for every $x\in \mX$,
then
$$
\lim_{j\to\be}\big\|\mI_k^{\mu,\nu_j}\big\|_{\mB_{\mX,\mY}^{\mu,\nu_j}}
=\big\|\mI_k^{\mu,\nu}\big\|_{\mB_{\mX,\mY}^{\mu,\nu}}.
$$
\end{lem}

\Pf Clearly, $\mI_k^{\mu,\nu_j}\big(\mI_k^{\mu,\nu_j}\big)^*$ is an integral operator on $L^2(\mX,\mu)$
with kernel $l_j(x,y)=\int_{\mY}k(x,t)\ov{k(y,t)}\,d\nu_j(t)$. Besides, the sequence
$\{l_j\}$ converges in $L^2(\mX\times\mX,\mu\otimes\mu)$ to the function $l$ defined by
$l(x,y)=\int_{\mY}k(x,t)\ov{k(y,t)}\,d\nu(t)$, which is the kernel of the integral operator
$\mI_k^{\mu,\nu}\big(\mI_k^{\mu,\nu}\big)^*$.
Hence,
\begin{align*}
\lim_{j\to\be}\big\|\mI_k^{\mu,\nu_j}\big\|_{\mB_{\mX,\mY}^{\mu,\nu_j}}^2
&=\lim_{j\to\be}\big\|\mI_k^{\mu,\nu_j}\big(\mI_k^{\mu,\nu_j}\big)^*\big\|_{
\mB_{\mX,\mY}^{\mu,\nu_j}}\\[.2cm]
&=\big\|\mI_k^{\mu,\nu}\big(\mI_k^{\mu,\nu}\big)^*\big\|_{\mB_{\mX,\mY}^{\mu,\nu}}
=\big\|\mI_k^{\mu,\nu}\big\|_{\mB_{\mX,\mY}^{\mu,\nu}}^2.\quad\bl
\end{align*}

\begin{cor}
\label{mlim}
Let $\mX$ and $\mY$ be compact subsets of $\R$, and
let $\mu$ and $\nu$ be finite positive Borel measures on $\mX$ and $\mY$.
Suppose that $\{\nu_j\}_{j=1}^\be$ is a sequence  of finite positive Borel measures on $\mY$ that
converges to $\nu$ in $\s\big((C(\mY))^*,C(\mY)\big)$.
If $\Phi$ is a Borel function on $\mX\times\mY$ such that
$\Phi(x,\cdot)\in C(\mY)$ for all $x\in \mX$,
then
$\|\Phi\|_{\frak M_{\mX,\mY}^{\,\mu,\nu}}\le\liminf\limits_{j\to\be}\|\Phi\|_{\frak M_{\mX,\mY}^{\,\mu,\nu_j}}$.
\end{cor}

\Pf It is easy to see that
$$
\|\Phi\|_{\frak M_{\mX,\mY}^{\,\mu,\nu}}
=\sup\big\{\|\Phi k\|_{\mB_{\mX,\mY}^{\,\mu,\nu}}:
~k\in C(\mX\times\mY),~\|k\|_{\mB_{\mX,\mY}^{\,\mu,\nu}}\le1\big\}.
$$
Let $k\in C(\mX\times\mY)$ with $\|k\|_{L^2(\mu\otimes\nu)}>0$. Then
\begin{align*}
\|\Phi k\|_{\mB_{\mX,\mY}^{\,\mu,\nu}}&=\lim_{j\to\be}\|\Phi k\|_{\mB_{\mX,\mY}^{\,\mu,\nu_j}}
\le\liminf_{j\to\be}\Big(\|\Phi\|_{\frak M_{\mX,\mY}^{\,\mu,\nu_j}}\|k\|_{\mB_{\mX,\mY}^{\,\mu,\nu_j}}\Big)\\[.2cm]
&=\liminf_{j\to\be}\|\Phi\|_{\frak M_{\mX,\mY}^{\,\mu,\nu_j}}\lim_{j\to\be}\|k\|_
{\mB_{\mX,\mY}^{\,\mu,\nu_j}}
=\|k\|_{\mB_{\mX,\mY}^{\,\mu,\nu}}\liminf_{j\to\be}\|\Phi\|_{\frak M_{\mX,\mY}^{\,\mu,\nu_j}}
\end{align*}
which implies the result. $\bl$

We are going to use the following notation: for a measure $\mu$ and an integrable function $\f$, we write $\nu=\f\mu$ if $\nu$ is the (complex) measure defined by $d\nu=\f\,d\mu$.

The following fact can be proved very easily.

\begin{lem}
\label{dapr}
Let $\nu$ and $\nu_0$ be finite Borel measures on $\R$ with compact
supports. Suppose that $\supp\nu_0\subset\supp\nu$. Then there exists a sequence
$\{\f_j\}_{j=1}^\be$ in $C(\R)$ such that $\f_j\ge0$ everywhere on $\R$ for all $j$
and $\nu_0=\lim\limits_{j\to\be}\f_j\nu$ in  $\s\big(\big(C(\supp\nu)\big)^*,C(\supp\nu)\big)$.
\end{lem}

\medskip

{\bf Proof of Theorem \ref{2mult+}.} Put $\mX_n\df[-n,n]\cap X$ and
$\mY_n\df[-n,n]\cap Y$. Clearly, $\left\{\|\Phi\|_{\frak M_{\mX_n,\mY_n}^{\,\mu,\nu}}\right\}$ is a nondecreasing sequence and
$$
\lim_{n\to\be}
\|\Phi\|_{\frak M_{\mX_n,\mY_n}^{\,\mu,\nu}}=
\|\Phi\|_{\frak M_{\mX,\mY}^{\,\mu,\nu}}.
$$
This allows us to reduce the general case
to the case when $\mX$ and $\mY$ are compact.
Besides, it suffices to consider the case where $\mu_0=\mu$.
Indeed, the case $\nu_0=\nu$ can be reduced to the case $\mu_0=\mu$, and we have
$$
\|\Phi\|_{\frak M_{\mX,\mY}^{\,\mu_0,\nu_0}}\le
\|\Phi\|_{\frak M_{\mX,\mY}^{\,\mu,\nu_0}}\le\|\Phi\|_{\frak M_{\mX,\mY}^{\,\mu,\nu}}.
$$
Let $\mX$ and $\mY$ be compact, and $\mu=\mu_0$. Applying Lemma \ref{dapr},
we can take a sequence $\{\f_j\}_{j=1}^\be$ of nonnegative functions in $C(\R)$ such that
$\nu_0=\lim\limits_{j\to\be}\f_j\nu$ in the weak topology $\s((C(\mY)^*,C(\mY))$.
Put $\nu_j\df\f_j\nu$.
By Theorem \ref{2mult}, $\|\Phi\|_{\frak M_{\mX,\mY}^{\,\mu,\nu_j}}
\le\|\Phi\|_{\frak M_{\mX,\mY}^{\,\mu,\nu}}$
for every $j\ge1$. It remains to apply Corollary \ref{mlim}. $\bl$

Theorem \ref{2mult+} implies the following fact:

\begin{thm}
\label{36}
Let $\mX$ and $\mY$ be closed subsets of $\R$ and let $\Phi$ be a
function on $\mX\times\mY$ that is continuous in each variables.
Suppose that $\mu$ and $\nu$ are positive regular Borel measures on $\mX$ and $\mY$
such that $\supp\mu=\mX$ and $\supp\nu=\mY$.
Then $\|\Phi\|_{\frak M_{\mX,\mY}}=\|\Phi\|_{\frak M_{\mX,\mY}^{\,\mu,\nu}}$.
\end{thm}

%

The following result is well known.

\medskip

{\it Let $f\in C(\R)$. Put $\Phi(x,y)\df f(x-y)$. Then
$\Phi\in\frak M_{\R,\R}$ if and only if $f$ is the Fourier transform of
a complex measure on $\R$. Moreover,
$\|\Phi\|_{\frak M_{\R,\R}}=|\mu|(\R)$.}

\medskip

A similar statement holds for any locally compact abelian group. In particular,
it is true for the group $\Z$:

{\it Let $f$ be a function defined on $\Z$. Put $\Phi(m,n)\df f(m-n)$. Then
$\Phi\in\frak M_{\Z,\Z}$ if and only if $\{f(n)\}_{n\in\Z}$ are
the Fourier coefficients of a complex Borel measure $\mu$
on the unit circle $\T$. Moreover,
$\|\Phi\|_{\frak M_{\Z,\Z}}=|\mu|(\T)$.}

\medskip

We need the following well-known fact.

\begin{lem}
\label{Hmn}
Let
$$
H(m,n)\df\left\{\begin{array}{ll}\frac1{m-n},&\text {if}\,\,\,\,m,n\in\Z,~ m\ne n,\\[.2cm]
0,&\text {if}\,\,\,\,m=n\in\Z.
\end{array}\right.
$$
Then $\|H\|_{\frak M_{\Z,\Z}}=\dfrac\pi2$.
\end{lem}

\Pf It suffices to observe that
$$
H(n,0)=\frac1{2\pi}\int_0^{2\pi}{\rm i}(\pi-t)e^{-{\rm i}nt}\,dt\quad\text{and}\quad
\frac1{2\pi}\int_0^{2\pi}|\pi-t|\,dt=\frac\pi2.\quad\bl
$$

\

\section{\bf Remarks on absolutely convergent Fourier integrals}
\setcounter{equation}{0}
\label{FouIn}

\

In this section we collect elementary estimates of certain functions in the space of absolutely convergent Fourier integrals. Such estimates will be used in the next section for estimates of certain functions in the space of Schur multipliers.

We are going to deal with the space
$$\widehat L^1=\widehat L^1(\R)\df\mathscr F(L^1(\R)),\quad
\|f\|_{\widehat L^1}=\|f\|_{\widehat L^1(\R)}\df\big\|\mathscr F^{-1} f\big\|_{L^1}.
$$
Here we use the notation $\F$ for {\it Fourier transform}:
$$
(\F f)(x)\df\int_\R f(t)e^{-{\rm i}xt}\,dt,\quad f\in L^1(\R).
$$

Unless otherwise stated, an interval throughout the paper means
a closed nondegenerate (not necessarily finite) interval.
For such an interval $J$, we consider the class $\widehat L^1(J)$ defined by
$\widehat L^1(J)\df\big\{f\big|J:f\in \widehat L^1\big\}$.
If $f\in C(J)$, we put
$$
\|\f\|_{\widehat L^1(J)}\df\inf\big\{\|f\|_{\widehat L^1}:f\big|J=\f\big\}.
$$
For $\f\in C(\R)$, we put $\|\f\|_{\widehat L^1(J)}\df\big\|\f\big|J\big\|_{\widehat L^1(J)}$.
Clearly, $\|\f\|_{L^\be(J)}\le\|\f\|_{\widehat L^1(J)}$.

For an interval $J$, {\it we use the notation $|J|$ for its length}.

It is easy to see that the nonzero constants belong to the space $\widehat L^1(J)$
for bounded intervals $J$ and $\|\mathbf1\|_{\widehat L^1(J)}=1$.
Moreover,
$$
\widehat L^1(J)=\big\{(\mathscr F\mu)|J:\mu\in \M(\R)\big\}\quad
\mbox{and}\quad \|f\|_{\widehat L^1(J)}=\inf\big\{\|\mu\|_{\M}:(\mathscr F\mu)\big|J=f\big\}
$$
for every bounded interval $J$,
where $\M(\R)$ denotes the space of (complex) Borel measures on $\R$.

In this section we are going to discuss (mostly known) estimates for $\|\cdot\|_{\widehat L^1(J)}$.

First, we recall the P\'olya theorem, see \cite{Po}.

\medskip

{\it Let $f$ be an even continuous functions such that $f\big|[0,\be)$ is decreasing convex function
vanishing at the infinity. Then $f\in\widehat L^1$ and $\|f\|_{\widehat L^1}=f(0)$.}

\medskip

This theorem readily implies the following fact.

\begin{lem}
\label{Polya}
Let $f$ be a continuous function on a closed ray $J$ that vanishes at infinity.
Suppose that $f$ is monotone and convex (or concave). Then $f\in\widehat L^1(J)$
and $\|f\|_{\widehat L^1(J)}=\max\limits_J|f|$.
\end{lem}

In what follows by a {\it locally absolutely continuous function} on $\R$ we mean a function whose restriction to any compact interval is absolutely continuous.

\begin{lem}
Let $f$ be a locally absolutely continuous function in $L^2(\R)$ such that $f^\prime\in L^2(\R)$. Then $f\in\widehat L^1(\R)$ and $\|f\|_{\widehat L^1}^2\le\|f\|_{L^2}\|f^\prime\|_{L^2}$.
\end{lem}

\Pf Put $a=\|f\|_{L^2}$, $b=\|f^\prime\|_{L^2}$.
By Plancherel's theorem,
$$
\big\|\F^{-1}f\big\|_{L^2}^2=\frac{a^2}{2\pi}\quad\text{and}\quad
\big\|x\F^{-1}f\big\|_{L^2}^2=\frac{b^2}{2\pi}.
$$
Hence,
$$
\Big\|\sqrt{b^2+a^2x^2}\,\,\mathscr F^{-1}f\Big\|_{L^2}^2=\frac{a^2b^2}{\pi}.
$$
and
by the Cauchy--Bunyakovsky inequality,
\bey
\big\|\mathscr F^{-1}f\big\|_{L^1}\le\frac{ab}{\sqrt{\pi}}\left\|\frac1{\sqrt{a^2x^2+b^2}}\right\|_{L^2}
=\sqrt{ab}.
\quad\bl
\eey

\begin{cor}
\label{Hxy}
Let $a>0$.
Put
$$
f_a(x)\df\left\{\begin{array}{ll}a^{-2}x,&\text {if}\,\,\,|x|\le a,\\[.2cm]
x^{-1},&\text {if}\,\,\,\,|x|\ge a.
\end{array}\right.
$$
Then $f_a\in\widehat L^1(\R)$ and $\|f_a\|_{\widehat L^1}\le\frac2a$.
\end{cor}

\Pf It suffices to observe that $\|f_a\|_{L^2}^2=\frac8{3a}$,
$\|f_a^\prime\|_{L^2}^2=\frac8{3a^3}$, and $\sqrt{\frac83}\le2$. $\bl$

\begin{lem}
\label{lip}
Let $J$ be a bounded interval and let $f$ be a Lipschitz function on $\R$ such that $\supp f\subset J$. Then
$f\in\widehat L^1$ and
$$
\|f\|_{\widehat L^1}\le\frac1{\root4 \of {12}}|J|\cdot\|f^\prime\|_{L^\be}.
$$
\end{lem}

\Pf Let $J=[-a,a]$. Clearly, $|f(x)|\le(a-|x|)\|f^\prime\|_{L^\be}$ for all
$x\in J$. Hence,
$$
\|f\|_{L^2}^2\le2\|f^\prime\|_{L^\be}^2\int_0^a(a-t)^2\,dt=\frac1{12}\|f^\prime\|_{L^\be}^2|J|^3.
$$
Using the obvious inequality $\|f^\prime\|_{L^2}^2\le\|f^\prime\|_{L^\be}^2|J|$,
we get the desired estimate. $\bl$

\begin{cor}
\label{13}
Let $f$ be a Lipschitz function on $\R$ such that $f(0)=0$.
Then \lb $\|f\|_{\widehat L^1(J)}\le\dfrac2{\root 4 \of {12}}|J|\cdot\|f^\prime\|_{L^\be}$
for every bounded interval J that contains $0$.
\end{cor}

\Pf Put $2J\df\{2x:~x\in J\}$.
Clearly, there exists a function $f_J$ in $C(\R)$ such that $f_J=f$ on $J$,
$\supp f_J\subset2J$, and $\|f_J^\prime\|_{L^\be}\le\|f^\prime\|_{L^\be}$. $\bl$

\begin{lem}
\label{logaa-}
Let $f$ be a locally absolutely continuous function on $\R$ such that
\lb$(1+|x|)f^\prime(x)\in L^2(\R)$. Suppose that $\lim\limits_{x\to-\be}f(x)=0$
and $\lim\limits_{x\to\be}f(x)=1$. Then
$$
\|f\|_{\widehat L^1(-\be,a]}\le\frac1{\sqrt{\pi}}\|f^\prime\|_{L^2}+\sqrt{\frac2\pi}\|xf^\prime\|_{L^2}
+\frac{7}{2\pi}+\frac2\pi\log a.
$$
for every $a\ge2$.
\end{lem}
\Pf Put
$$
f_a(x)\df f(x)-a^{-1}\int_{-\be}^x\chi_{[a,2a]}(t)\,dt.
$$
Clearly, $\|f\|_{\widehat L^1(-\be,a]}\le\|f_a\|_{\widehat L^1}$.

We have
$$
-{\rm i}x\mathscr F^{-1}f_a=\mathscr F^{-1}(f_a^\prime)=\mathscr F^{-1}
(f^\prime)-\frac{e^{2a{\rm i}x}-e^{a{\rm i}x}}{2\pi a{\rm i}x}
$$
Put $h\df\mathscr F^{-1}(f^\prime)$. Then
\begin{align*}
\|f_a\|_{\widehat L^1}=&
\int_\R\left|h(x)-\frac{e^{2a{\rm i}x}-e^{a{\rm i}x}}{2\pi a{\rm i}x}\right|\cdot\frac{dx}{|x|}\\[.2cm]
\le&\int_{-1}^1\frac{|h(x)-h(0)|}{|x|}\,dx
+\frac1{2\pi}\int_{-1}^1\left|\frac{e^{2a{\rm i}x}-e^{a{\rm i}x}}{a{\rm i}x}-1\right|\cdot\frac{dx}{|x|}\\[.2cm]
&+\int_{\{|x|\ge1\}}\frac{|h(x)|}{|x|}\,dx
+\frac1{2\pi a}\int_{\{|x|\ge1\}}\frac{|e^{a{\rm i}x}-1|}{x^2}\,dx.
\end{align*}

We have
$$
\int_0^1\frac{|h(x)-h(0)|}{x}\,dx\le\int_0^1\frac1x\left(\int_0^x|h^\prime(t)|\,dt\right)dx=
\int_0^1|h^\prime(t)|\cdot|\log t|\,dt.
$$
Hence,
\begin{align*}
\int_{-1}^1\frac{|h(x)-h(0)|}{|x|}\,dx&\le\int_{-1}^1|h^\prime(t)|\cdot\big|\log |t|\,\big|\,dt\\[.2cm]
&\le\|h^\prime\|_{L^2}\left(\int_{-1}^1\log^2|t|\,dt\right)^{1/2}=\sqrt{\frac2\pi}\,\,\|xf^\prime(x)\|_{L^2}
\end{align*}
because $h^\prime=\mathscr F^{-1}({\rm i}xf^\prime)$.

By Taylor's formula for the function $e^{2{\rm i}x}-e^{{\rm i}x}$, we have
$$
\big|e^{2{\rm i}x}-e^{{\rm i}x}-{\rm i}x\big|\le\frac52\,{x^2}.
$$
Thus
\begin{align*}
\frac1{2\pi}\int_{-1}^1\left|\frac{e^{2a{\rm i}x}-e^{a{\rm i}x}}{a{\rm i}x}-1\right|\cdot\frac{dx}{|x|}
&=\frac1{2\pi}\int_{-a}^a\left|\frac{e^{2{\rm i}x}-e^{{\rm i}x}}{{\rm i}x}-1\right|\cdot\frac{dx}{|x|}\\[.2cm]
&\le\frac1{2\pi}\int_{-a}^a\min\left\{\frac52,\frac2{|x|}\right\}\,dx
\le\frac1{2\pi}(5+4\log a).
\end{align*}
Finally,
$$
\int_{|x|\ge1}\frac{|h(x)|}{|x|}\,dx\le\sqrt2\|h\|_{L^2}=\frac1{\sqrt{\pi}}\,\,\|f^\prime\|_{L^2}
$$
by the Cauchy--Bunyakovsky inequality and
$$
\frac1{2\pi a}\int_{\{|x|\ge1\}}\frac{|e^{a{\rm i}x}-1|}{x^2}\,dx
=\frac1{2\pi}\int_{\{|x|\ge a\}}\frac{|e^{{\rm i}x}-1|}{x^2}\,dx\le\frac2{\pi a}\le\frac1\pi
$$
for $a\ge2$.
This implies the desired inequality. $\bl$

\begin{thm}
\label{smalJ}
Let $J$ be a bounded interval containing $0$. Then
\bay
\label{smJ}
\left\|\frac{e^x-1}{e^x+1}\right\|_{\widehat L^1(J)}\le\frac1{\root 4 \of{12}}|J|\le\frac35|J|.
\ey
\end{thm}

\Pf Ii suffices to observe that $\big\|(\frac{e^x-1}{e^x+1})^\prime\big\|_{L^\be}=\frac12$
and apply Corollary \ref{13}. $\bl$

Theorem \ref{smalJ} gives a sufficiently sharp estimate of $\widehat L^1$-norm for little intervals $J$.
For big intervals $J$, this estimate will be improved in Corollary \ref{bigJ}.

\begin{thm}
\label{35}
Let $a\ge2$. Then
$$
\left\|\frac{e^x}{1+e^x}\right\|_{\widehat L^1(-\be,a]}\le2+\frac2\pi\log a.
$$
\end{thm}

\Pf
We have
$$
\left\|\left(\frac{e^x}{e^x+1}\right)^\prime\right\|_{L^2}^2=\int_\R\frac{e^{2x}\,dx}{(e^x+1)^4}=
\int_0^\be\frac{t\,dt}{(t+1)^4}=\frac16,
$$
and\footnote{In fact, $\|x(\frac{e^x}{e^x+1})^\prime\|_{L^2}^2=\frac{\pi^2}{18}-\frac13$.}
$$
\left\|x\left(\frac{e^x}{e^x+1}\right)^\prime\right\|_{L^2}^2=2\int_0^\be\frac{x^2e^{2x}\,dx}{(e^x+1)^4}\le
2\int_0^\be x^2e^{-2x}\,dx=\frac12,
$$
whence for $a\ge2$,
$$
\left\|\frac{e^x}{1+e^x}\right\|_{\widehat L^1(-\be,a]}
\le\frac1{\sqrt{6\pi}}+\frac1{\sqrt{\pi}}
+\frac{7}{2\pi}+\frac2\pi\log a
\le2+\frac2\pi\log a
$$
by Lemma \ref{logaa-}. $\bl$

\medskip

{\bf Remark.} Lemma \ref{Polya} implies that
$$
\left\|\frac{e^x}{1+e^x}\right\|_{\widehat L^1(-\be,a]}\le\frac{e^a}{1+e^a}\le e^a
$$
for $a\le0$ but we do not need this inequality.

\medskip

\begin{cor}
\label{bigJ}
Let $J$ be a bounded interval containing $0$. Then
$$
\left\|\frac{e^x-1}{e^x+1}\right\|_{\widehat L^1(J)}\le5+\frac4\pi\log\left(\frac12|J|\right)
$$
if $|J|\ge4$.
\end{cor}
\Pf We may assume that the center of $J$ is nonpositive. Then
$J\subset\big(-\be,\frac12|J|\big]$. We have
$$
\left\|\frac{e^x-1}{e^x+1}\right\|_{\widehat L^1(J)}
\le1+2\left\|\frac{e^x}{e^x+1}\right\|_{\widehat L^1(J)}
\le5+\frac4\pi\log a
=5+\frac4\pi\log\left(\frac12|J|\right).\quad\bl
$$

\

\section{\bf Estimates of certain multiplier norms}
\label{Estim}

\

In this section we are going to obtain sharp estimates
for the Schur multiplier norms
\bay
\label{tens}
\left\|\frac{e^x-e^y}{e^x+e^y}\right\|_{\frak M_{J_1,J_2}}
=\left\|\frac{e^{x-y}-1}{e^{x-y}+1}\right\|_{\frak M_{J_1,J_2}}
\ey
for every intervals $J_1$ and $J_2$. First, we consider two special cases. In the first case $J_1=J_2$ while in the second case $J_1$ and $J_2$ do not overlap, i.e., their intersection has at most one point.

\begin{thm}
\label{Jdis}
Let $J_1$ and $J_2$ be nonoverlapping intervals. Then
$$
\left\|\frac{e^x-e^y}{e^x+e^y}\right\|_{\frak M_{J_1,J_2}}\le2.
$$
\end{thm}

\Pf Clearly, either $J_1-J_2\subset(-\be,0]$ or $J_1-J_2\subset[0,\be)$.
It suffices to consider the case when $J_1-J_2\subset(-\be,0]$.
Then
\begin{align*}
\left\|\frac{e^x-e^y}{e^x+e^y}\right\|_{\frak M_{J_1,J_2}}
&\le1+2\left\|\frac{e^{x-y}}{e^{x-y}+1}\right\|_{\frak M_{J_1,J_2}}\\[.2cm]
&\le1+2\left\|\frac{e^{x}}{e^{x}+1}\right\|_{\widehat L^1(-\be,0]}
=2
\end{align*}
by the P\'olya theorem \cite{Po}, see also Lemma \ref{Polya}. $\bl$

\begin{thm}
\label{Jogr}
Let $J$ be a bounded interval.
Then
$$
\left\|\frac{e^x-e^y}{e^x+e^y}\right\|_{\frak M_{J,J}}\le
\min\left\{\frac65|J|\,,\,\,\,\,5+\frac4\pi\log_+|J|\right\}
$$
and so
$$
\left\|\frac{e^x-e^y}{e^x+e^y}\right\|_{\frak M_{J,J}}\le4\log(1+|J|).
$$
\end{thm}

\Pf We have
$$
\left\|\frac{e^x-e^y}{e^x+e^y}\right\|_{\frak M_{J,J}}
\le\left\|\frac{e^x-1}{e^x+1}\right\|_{\widehat L^1(J-J)}.
$$
Note that $|J-J|=2|J|$ and $0\in J-J$. The result follows now from
Theorem \ref{smalJ} and Corollary \ref{bigJ}. $\bl$

\begin{thm}
\label{Jdis+}
Let $J_1$ and $J_2$ be nonoverlapping intervals and
let $J$ be the convex hull of $J_1\cup J_2$.
Then
$$
\frac{e-1}{e+1}\min\big\{1,|J|\big\}\le\left\|\frac{e^x-e^y}{e^x+e^y}\right\|_{\frak M_{J_1,J_2}}\le
\min\left\{2,\frac65|J|\right\}.
$$
\end{thm}

\Pf The upper estimate follows readily from Theorems \ref{Jdis} and \ref{Jogr}.
Let us prove the lower estimate. We have
\bey
\left\|\frac{e^x-e^y}{e^x+e^y}\right\|_{\frak M_{J_1,J_2}}
\ge
\sup_{x\in J_1,\,y\in J_2}\left|\frac{e^x-e^y}{e^x+e^y}\right|\ge
\frac{e^{|J|}-1}{e^{|J|}+1}\ge\frac{e-1}{e+1}\min\{1,|J|\}
\eey
because the function $t\mapsto\frac{e^t-1}{t(e^t+1)}$ decreases on $[0,\be)$, while the function $t\mapsto\frac{e^t-1}{e^t+1}$ increases.
$\bl$

\begin{thm}
\label{Jogr-}
Let $J$ be a bounded interval.
Then
$$
\left\|\frac{e^x-e^y}{e^x+e^y}\right\|_{\frak M_{J,J}}\ge
\frac19\min\Big\{|J|\,,\,\,\,1+\log_+|J|\Big\}.
$$
\end{thm}

\Pf Put $Q_\e(t)\df\frac1\pi\frac{t}{t^2+\e^2}$, where $\e>0$. Let us consider
the convolution operator $\mathcal C_{Q_\e}$ on $L^2(\R)$, $\mathcal C_{Q_\e}f\df f*Q_\e$.
Clearly, $\|\mathcal C_{Q_\e}\|=\|\mathscr F Q_\e\|_{L^\be}=1$, see,
for example, \cite{Ga}, Ch. III, \S\,1.
Note that $\mathcal C_{Q_\e}$ is an integral operator with kernel $Q_\e(x-y)$.
We can define the integral operator $X_{J,\e}$ on $L^2(J)$ with
kernel
$$
\frac1\pi\,\frac{x-y}{(x-y)^2+\e^2}\,\,\frac{e^x-e^y}{e^x+e^y}.
$$
We have
\begin{align*}
|J|\cdot\|X_{J,\e}\|\ge(X_{J,\e}\chi_J,\chi_J)&=\frac1\pi\iint\limits_{J\times J}\frac{x-y}{(x-y)^2+\e^2}\,\,
\frac{e^x-e^y}{e^x+e^y}\,dxdy\\[.2cm]
&=\frac{2\sqrt2}\pi\int_0^{|J|}\frac{t}{t^2+\e^2}\,\,
\frac{e^t-1}{e^t+1}(|J|-t)\,dt
\end{align*}
and
$$
\|X_{J,\e}\|\le\|\mathcal C_{Q_\e}\|\cdot\left\|\frac{e^x-e^y}{e^x+e^y}\right\|_{\frak M_{J,J}}
=\left\|\frac{e^x-e^y}{e^x+e^y}\right\|_{\frak M_{J,J}}.
$$
Hence,
$$
\left\|\frac{e^x-e^y}{e^x+e^y}\right\|_{\frak M_{J,J}}
\ge\frac{2\sqrt2}\pi\cdot\frac{1}{|J|}\int_0^{|J|}\frac{t}{t^2+\e^2}\,\,
\frac{e^t-1}{e^t+1}(|J|-t)\,dt
$$
for every $\e>0$, whence
$$
\left\|\frac{e^x-e^y}{e^x+e^y}\right\|_{\frak M_{J,J}}
\ge\frac{2\sqrt2}\pi\int_0^{|J|}\frac{e^t-1}{t(e^t+1)}\left(1-\frac{t}{|J|}\right)\,dt
\ge\frac{\sqrt2}\pi\int_0^{|J|}\frac{e^t-1}{t(e^t+1)}\,dt
$$
because the function $t\mapsto\frac{e^t-1}{t(e^t+1)}$ decreases on $(0,\be)$. It follows that
$$
\left\|\frac{e^x-e^y}{e^x+e^y}\right\|_{\frak M_{J,J}}
\ge\frac{\sqrt2}\pi\cdot\frac{e-1}{e+1}\int_0^{|J|}\min\{1,t^{-1}\}\,dt.
$$
This implies the desired estimate. $\bl$

\medskip

{\bf Remark 1.} Every rectangle $J_1\times J_2$ is the union at most of
three rectangles, each of which satisfies the hypotheses of either Theorem \ref{Jogr}
or Theorem \ref{Jdis+}. This allows us to obtain a sharp estimate for the norms in \rf{tens}
for every rectangle $J_1\times J_2$.

\medskip

{\bf Remark 2.} Remark 1 and the change of variables $x\mapsto\log x$, $y\mapsto\log y$
allow us to obtain a sharp estimate for $\left\|\dfrac{x-y}{x+y}\right\|_{\frak M_{J_1,J_2}}$,
where $J_1$ and $J_2$ are intervals containing in $(0,\be)$.

\medskip

We proceed now to estimates of multiplier norms that will be used in this paper.

\begin{thm}
\label{0albe}
There exists a positive number $C$ such that
$$
\left\|\frac{e^x-e^y}{e^x+e^y}\right\|_{\frak M_{[a,\be),(-\be,b]}}\le C\log(2+(b-a)_+)
$$
for all $a,b\in\R$.
\end{thm}

\Pf The result follows from Theorems \ref{Jdis} if $a\ge b$.
If $a<b$, then
$$
[a,\be)\times(-\be,b]=
([a,b]\times[a,b])\cup([a,b]\times(-\be,a])\cup([b,\be)\times(-\be,b]),
$$
and we can apply Theorem \ref{Jogr}  to the first rectangle and
Theorem \ref{Jdis} to the remaining rectangles. $\bl$

\begin{thm}
\label{0albe+}
There exists a positive number $C$ such that
$$
\left\|\frac{e^x-e^y}{e^x+e^y}\right\|_{\frak M_{\R,[a,b]}}\le C\log(2+b-a)
$$
for every $a,b\in\R$ satisfying $a<b$.
\end{thm}

\Pf We have
$$
\R\times[a,b]=([a,b]\times[a,b])\cup((-\be,a]\times[a,b])\cup([b,\be)\times[a,b]).
$$
It remains to apply Theorem \ref{Jogr}  to the first rectangle and
Theorem  \ref{Jdis} to the remaining rectangles. $\bl$

%
%

\begin{thm}
\label{albe}
There exists a positive number $c$ such that
$$
\left\|\frac{x-y}{x+y}\right\|_{\frak M_{[a,\be),[0,b]}}\le c\log\left(2+\log_+\frac ba\right)
$$
for all $a,b\in(0,\be)$.
\end{thm}

\Pf Theorem \ref{0albe} with the help of
the change of variables $x\mapsto\log x$ and $y\mapsto\log y$
yields
$$
\left\|\frac{x-y}{x+y}\right\|_{\frak M_{[a,\be),[\e,b+\e]}}\le c\log\left(2+\log_+\frac{b+\e}a\right)
$$
for every $\e>0$, whence
$$
\left\|\frac{x-y-\e}{x+y+\e}\right\|_{\frak M_{[a,\be),[0,b]}}\le c\log\left(2+\log_+\frac{b+\e}a\right)
$$
for every $\e>0$. It remains to pass to the limit as $\e\to0$. $\bl$

\begin{thm}
\label{albe+}
There exists a positive number $c$ such that
$$
\left\|\frac{x-y}{x+y}\right\|_{\frak M_{[a,b],[0,\be)}}\le c\log\left(2+\log\frac ba\right)
$$
whenever $a,b\in(0,\be)$ and $a<b$.
\end{thm}

\Pf The result follows from Theorem \ref{0albe+} in the same way as
Theorem \ref{albe} follows from Theorem \ref{0albe}. $\bl$

\begin{thm}
\label{albe-}
There exists a positive number $c$ such that
$$
\left\|\frac{x-y}{x+y}\right\|_{\frak M_{[a,b],[a,b]}}\ge c\log\left(1+\log\frac ba\right)
$$
whenever $a,b\in(0,\be)$ and $a<b$.
\end{thm}

\Pf The result follows from Theorem \ref{Jogr-}  with the help of
the change of variables $x\mapsto\log x$ and $y\mapsto\log y$. $\bl$

\
\section{\bf Operator Lipschitz functions and operator modulus of continuity}
\setcounter{equation}{0}
\label{OL}

\

In this section we study operator Lipschitz functions on closed subsets of the real line. It is well known that a function $f$ on $\R$ is operator Lipschitz if and only if it is {\it commutator Lipschitz}, i.e.,
$$
\|f(A)R-Rf(A)\|\le\const\|AR-RA\|
$$
for an arbitrary bounded operator $R$ and an arbitrary self-adjoint operator $A$.

It turns out that the same is true for functions on closed subsets of $\R$;
moreover the operator Lipschitz constant coincides with the commutator Lipschitz constant.
The following theorem was proved in \cite{AP2} (Th. 10.1)
in the case $\frak F=\R$.
The general case is analogous to the case $\frak F=\R$.
See also \cite{KS} where similar results for symmetric
ideal norms are considered.

\begin{thm}
\label{Olip}
Let $f$ be a continuous function defined on a closed subset $\frak F$ of $\R$ and let $C\ge0$. The following are equivalent:

{\em(i)} $\|f(A)-f(B)\|\le C\|A-B\|$ for arbitrary self-adjoint operators $A$ and $B$
with spectra in $\frak F$;

{\em(ii)} $\|f(A)R-Rf(A)\|\le C\|AR-RA\|$ for all self-adjoint operators $A$ with $\s(A)\subset \frak F$ and all bounded operators $R$;

{\em(iii)} $\|f(A)R-Rf(B)\|\le C\|AR-RB\|$ for arbitrary self-adjoint operators $A$ and $B$
with spectra in $\frak F$
and for an arbitrary bounded operator $R$.
\end{thm}

A function $f\in C(\frak F)$ is said to be {\it operator
Lipschitz} if it satisfies the equivalent statements of Theorem \ref{Olip}.
We denote the set of operator Lipschitz functions on $\frak F$ by ${\rm OL}(\frak F)$.
For $f\in{\rm OL}(\frak F)$, we define $\|f\|_{{\rm OL}(\frak F)}$ to be the smallest constant satisfying the equivalent statements of Theorem \ref{Olip}.
Put $\|f\|_{{\rm OL}(\frak F)}=\be$ if $f\not\in {\rm OL}(\frak F)$.

It is well known that every $f$ in ${\rm OL}(\frak F)$ is differentiable at every non-isolated point of $\frak F$, see \cite{JW}.
Moreover, the same argument gives differentiability at $\infty$
in the following sense: there exists a finite limit $\lim\limits_{|x|\to+\infty}x^{-1}f(x)$
provided $\frak F$ is unbounded.

Let $f\in{\rm OL}(\frak F)$. Suppose that $\frak F$ has no isolated points.
Put
$$
\big(\dg f\big)(x,y)\df
\left\{\begin{array}{ll}\frac{f(x)-f(y)}{x-y},&\text {if}\,\,\,\,x,¸áy\in \frak F, \,\,\,x\ne y,\\[.2cm]
f^\prime(x),&\text {if}\,\,\,\,x\in \frak F, \,\,\,x=y.
\end{array}\right.
$$
The following equality holds:
\bay
\label{ravvo}
\|f\|_{{\rm OL}(\fF)}=\|\dg f\|_{{\fM}_{\fF,\fF}}.
\ey
The inequality $\|f\|_{{\rm OL}(\fF)}\le\|\dg f\|_{{\fM}_{\fF,\fF}}$ is an immediate consequence of the formula
\bay
\label{doi}
f(A)-f(B)=\iint(\dg f\big)(x,y)\,dE_A(x)(A-B)\,dE_B(y),
\ey
where $A$ and $B$ are self-adjoint operators with bounded $A-B$ whose spectra are in $\fF$, and $E_A$ and $E_B$ are the spectral measures of $A$ and $B$. The expression on the right is called a double operator integral. We refer the reader to \cite{BS1}, \cite{BS2}, and \cite{BS3} for the theory of double operator integrals elaborated by Birman and Solomyak. The validity of formula \rf{doi} under the assumption
$\dg f\in{\fM}_{\fF,\fF}$ and the inequality
$$
\left\|\iint(\dg f\big)(x,y)\,dE_A(x)(A-B)\,dE_B(y)\right\|\le\|\dg\|_{{\fM}_{\fF,\fF}}\|A-B\|
$$
was proved in \cite{BS3}. The opposite inequality in \rf{ravvo} is going to be proved in Corollary \ref{obr}.

In the general case for $f\in{\rm OL}(\frak F)$ we can define the function
$$
\big(\dg_0f\big)(x,y)\df
\left\{\begin{array}{ll}\frac{f(x)-f(y)}{x-y},&\text {if}\,\,\,\,x,y\in \frak F, \,\,\,x\ne y,\\[.2cm]
0,&\text {if}\,\,\,\,x\in \frak F, \,\,\,x=y.
\end{array}\right.
$$
The following inequalities hold:
\bay
\label{nva}
\|f\|_{{\rm OL}(\frak F)}\le\|\dg_0 f\|_{{\frak M}_{\frak F,\frak F}}\le
2\|f\|_{{\rm OL}(\frak F)}.
\ey
The first inequality in \rf{nva} follows from the formula
\bay
\label{doie}
f(A)-f(B)=\iint(\dg_0f\big)(x,y)\,dE_A(x)(A-B)\,dE_B(y),
\ey
whose validity can be verified in the same way as the validity of \rf{doi}. The second inequality in \rf{nva} is going to be verified in Corollary \ref{2OL}.

Let $f$ be a continuous function on a closed set $\frak F$, $\frak F\subset\R$.
We define the operator modulus of continuity $\O_{f, \frak F}$ as follows
$$
\O_{f,\frak F}(\d)\df\sup\big\{\|f(A)-f(B)\|:~A=A^*,~B=B^*,~\s(A),\s(B)\subset \frak F,~\|A-B\|\le\d\big\},
$$
and the commutator modulus of continuity as follows
$$
\O_{f,\frak F}^{\flat}(\d)\df\sup\big\{\|f(A)R-Rf(A)\|:~A=A^*,~\s(A)\subset \frak F,~\|R\|\le1,~\|AR-RA\|\le\d\big\}.
$$
One can prove that we get the same right-hand side if we require in addition that $R$ is self-adjoint.
On the other hand, $\|f(A)R-Rf(B)\|\le\O_{f,\frak F}^{\flat}\big(\|AR-RB\|\big)$ for every self-adjoint operators
with $\s(A),\s(B)\subset \frak F$ and for every bounded operator $R$ with $\|R\|\le1$.
Also, $\O_{f,\frak F}\le\O_{f,\frak F}^{\flat}\le2\O_{f,\frak F}$.

These results were obtained in \cite{AP2} in the case $\frak F=\R$. The same reasoning works in the general case.

\begin{lem}
\label{BOL}
Let $\frak F$ be a closed subset of $\R$ and
let $\mu$ and $\nu$ be regular positive Borel measures on $\frak F$.
Suppose that $k$ is a function in $L^2(\frak F\times\frak F,\mu\otimes\nu)$ such that $k=0$ on the diagonal $\D_\frak F\df\{(x,x):\,\,x\in\frak F\}$ almost everywhere  with respect to $\mu\otimes\nu$. Then
$$
\|k\,\dg_0 f\|_{\mB_{\fF,\fF}^{\,\mu,\nu}}\le\|f\|_{{\rm OL}(\frak F)}\|k\|_{\mB_{\fF,\fF}^{\,\mu,\nu}}
$$
for every continuous function $f$ on $\frak F$.
\end{lem}

\Pf Let $\fF_n\df\fF\cap[-n,n]$, and let $\mu_n$ and $\nu_n$ be the restrictions of $\mu$ and $\nu$ to $\fF_n$. Clearly,
$$
\lim\limits_{n\to\be}\|k\|_{{\mB}_{\fF_n,\fF_n}^{\mu_n,\nu_n}}
=\|k\|_{{\mB}_{\fF,\fF}^{\mu,\nu}}
\quad\mbox{for every}\quad k\in L^2(\fF\times\fF,\mu\otimes\nu)
$$
and
$$
\lim\limits_{n\to\be}\|f\|_{{\rm OL}(\fF_n)}=\|f\|_{{\rm OL}(\fF)}
\quad\mbox{for every}\quad f\in C(\fF).
$$
Thus we may assume that $\fF$ is compact.
It suffices to consider the case when $k$ vanishes in a neighborhood of
the diagonal $\D_\frak F$. Put $l(x,y)\df(x-y)^{-1}k(x,y)$.
Denote by $A$ and $B$ multiplications by the independent variable on $L^2(\frak F,\mu)$  and $L^2(\fF,\nu)$. Then
$\mI_k^{\mu,\nu}=A\mI_l^{\mu,\nu}-\mI_l^{\mu,\nu} B$ and
$\mI_{k\dg_0 f}^{\mu,\nu}=f(A)\mI_l^{\mu,\nu}-\mI_l^{\mu,\nu} f(B)$.
It remains to observe that
$$
\big\|f(A)\mI_l^{\mu,\nu}-\mI_l^{\mu,\nu} f(B)\big\|
\le\|f\|_{{\rm OL}(\fF)}\big\|A\mI_l^{\mu,\nu}-\mI_l^{\mu,\nu} B\big\|,
$$
$$
\big\|A\mI_l^{\mu,\nu}-\mI_l^{\mu,\nu} B\big\|=\|k\|_{\mB_{\fF,\fF}^{\,\mu,\nu}},
$$
and
$$
\big\|f(A)\mI_l^{\mu,\nu}-\mI_l^{\mu,\nu} f(B)\big\|=\|k\,\dg_0 f\|_{\mB_{\fF,\fF}^{\,\mu,\nu}}.
\quad\bl
$$

\begin{cor}
\label{cFF}
Let $\frak F$ be a closed subset of $\R$ with no isolated points, and
let $\mu$ and $\nu$ be finite positive Borel measures on $\fF$.
Suppose that $f$ is a differentiable function on $\frak F$ and $k\in L^2(\fF\times\fF,\mu\otimes\nu)$. If $k$ vanishes $\mu\otimes\nu$-almost everywhere
on the diagonal $\D_\frak F\df\{(x,x):\,\,x\in\frak F\}$, then
$$
\|k\,\dg f\|_{\mB_{\fF,\fF}^{\,\mu,\nu}}\le\|f\|_{{\rm OL}(\frak F)}\|k\|_{\mB_{\fF,\fF}^{\,\mu,\nu}}.
$$
\end{cor}

\Pf It suffices to observe that $k\,\dg f=k\,\dg_0 f$ almost everywhere
with respect to $\mu\otimes\nu$. $\bl$

\begin{cor}
\label{obr}
Let $\frak F$ be a closed subset of $\R$ with no isolated points, and
let $\mu$ and $\nu$ be finite positive Borel measures on $\frak F$.
If $f$ is a differentiable function on $\fF$, then
$$
\|\dg f\|_{{\frak M}_{\frak F,\frak F}}\le\|f\|_{{\rm OL}(\frak F)}.
$$
\end{cor}

\Pf Let $\mu$ be a regular Borel measure on $\frak F$ with no atoms
and such that $\supp\mu=\fF$. Then $(\mu\otimes\mu)(\D_\fF)=0$ and Corollary \ref{cFF}
implies that
$$
\|k\,\dg f\|_{\mB_{\fF,\fF}^{\,\mu,\mu}}\le\|f\|_{{\rm OL}(\frak F)}\|k\|_{\mB_{\fF,\fF}^{\,\mu,\mu}}
$$
for all $k\in L^2(\frak F\times\frak F,\mu\otimes\mu)$.
Hence, $\|\dg f\|_{\fM_{\fF,\fF}^{\,\mu,\mu}}\le\|f\|_{{\rm OL}(\frak F)}$.
It remains to apply Theorem \ref{36}. $\bl$

\begin{cor}
\label{2OL}
Let $\frak F$ be a closed subset of $\R$. Then
$$
\|\dg_0 f\|_{{\frak M}_{\frak F,\frak F}}\le2\|f\|_{{\rm OL}(\frak F)}
$$
for every $f\in C(\frak F)$.
\end{cor}

\Pf  Let $\mu$ and $\nu$ be
regular Borel measures on $\frak F$.
We have to verify that
$$
\|k\,\dg_0 f\|_{\mB_{\fF,\fF}^{\,\mu,\nu}}\le2\|f\|_{{\rm OL}(\fF)}\|k\|_{\mB_{\fF,\fF}^{\,\mu,\nu}}
$$
for every $k\in L^2(\frak F\times\frak F,\mu\otimes\nu)$.
Put
$k_0\df\chi_{\D_\frak F}k$ and $k_1\df k-k_0$. We have
$$
\|k_0\|_{\mB_{\fF,\fF}^{\,\mu,\nu}}\le\|k\|_{\mB_{\fF,\fF}^{\,\mu,\nu}}.
$$
This inequality can be verified easily. We leave the verification to the reader.

It follows that
$\|k_1\|_{\mB_{\fF,\fF}^{\,\mu,\nu}}\le\|k_0\|_{\mB_{\fF,\fF}^{\,\mu,\nu}}+\|k\|_{\mB_{\fF,\fF}^{\,\mu,\nu}}
\le2\|k\|_{\mB_{\fF,\fF}^{\,\mu,\nu}}$.
It remains to observe that
$$
\|k\,\dg_0 f\|_{\mB_{\fF,\fF}^{\,\mu,\nu}}=\|k_1\,\dg_0 f\|_{\mB_{\fF,\fF}^{\,\mu,\nu}}
\le\|f\|_{{\rm OL}(\fF)}\|k_1\|_{\mB_{\fF,\fF}^{\,\mu,\nu}}\le
2\|f\|_{{\rm OL}(\fF)}\|k\|_{\mB_{\fF,\fF}^{\,\mu,\nu}}.\quad\bl
$$

Let $\fF_1$ and $\fF_2$ be closed subsets of $\R$.
We define the space ${\rm OL}(\fF_1,\fF_2)$ as the space of
functions $f$ in $C(\fF_1\cup\fF_2)$ such that
\bay
\label{lip1R2}
\|f(A)R-Rf(B)\|\le C\|AR-RB\|
\ey
for all bounded operator $R$ and all self-adjoint operators $A$ and $B$ such that $\s(A)\subset\fF_1$
and $\s(B)\subset\frak F_2$ with some positive number $C$. Denote by $\|f\|_{{\rm OL}(\fF_1,\fF_2)}$
the minimal constant satisfying \rf{lip1R2}. Clearly,
$\|f\|_{{\rm OL}(\fF_1,\fF_2)}=\|f\|_{{\rm OL}(\fF_2,\fF_1)}$ and
$\|f\|_{{\rm OL}(\fF,\fF)}=\|f\|_{{\rm OL}(\fF)}$.
As in the case $\fF_1=\fF_2$, we can prove that
\bay
\label{nva+}
\|f\|_{{\rm OL}(\fF_1,\fF_2)}\le\|\dg_0 f\|_{{\frak M}_{\fF_1,\fF_2}}\le
2\|f\|_{{\rm OL}(\fF_1,\fF_2)}.
\ey
(cf. \rf{nva}).

\medskip

{\bf Remark.} In the case where $\fF_1\ne\fF_2$ we cannot claim that the inequality
\bay
\label{lip12}
\|f(A)-f(B)\|\le C\|A-B\|
\ey
for all self-adjoint $A$ and $B$  such that $\s(A)\subset\fF_1$
and $\s(B)\subset\frak F_2$ implies \rf{lip1R2}.

Indeed, in the case $f(t)=|t|$, $\fF_1=(-\be,0]$, and $\fF_2=[0,\be)$,
inequality \rf{lip12} holds with $C=1$ because
$$
\|A-B\|\le\|A+B\|
$$
for positive self-adjoin operators $A$ and $B$. However, inequality \rf{lip1R2} does not hold with any positive $C$. Indeed,
$$
\left\|\frac{|x|-|y|}{x-y}\right\|_{\fM_{(-\be,1],[1,\be)}}=
\left\|\frac{x-y}{x+y}\right\|_{\fM_{[1.\be),[1,\be)}}=\be
$$
by Theorem \ref{albe-}.

\medskip

\begin{thm}
\label{spectra}
Suppose that inequality {\em\rf{lip1R2}} holds for every bounded operator $R$
and arbitrary self-adjoint operators $A$ and $B$ with simple spectra such that
$\s(A)\subset\fF_1$ and $\s(B)\subset\frak F_2$. Then $f\in{\rm OL}(\fF_1,\fF_2)$
and $\|f\|_{{\rm OL}(\fF_1,\fF_2)}\le C$.
\end{thm}

\Pf We have to prove inequality \rf{lip1R2} for arbitrary self-adjoint operators $A$ and $B$ with $\s(A)\subset\fF_1$
and $\s(B)\subset\frak F_2$. It is convenient to think that the operators $A$ and $B$ act
in different Hilbert spaces. Let $A$ act in  $\mathcal H_1$ and
$B$ in $\mathcal H_2$. Then $R$ acts from $\mathcal H_2$
into $\mathcal H_1$. We are going to verify
that
$$
\big|(f(A)Ru,v)-(Rf(B)u,v)\big|=\big|(Ru,\ov f(A)v)-(f(B)u,R^*v)\big|\le C\|AR-RB\|
$$
for every unit vectors $u\in\mathcal H_2$ and $v\in\mathcal H_1$.
Denote by $\mathcal H_1^0$ and $\mathcal H_2^0$ the invariant subspaces of $A$ and $B$ generated by $v$ and $u$. Clearly, $A_0\df A|\mathcal H_1^0$ and $B_0\df B|\mathcal H_2^0$
are self-adjoint operators with simple spectra. Consider the operator
$R_0:\mathcal H_2^0\to\mathcal H_1^0$, $R_0h\df PRh$ for $h\in\mathcal H_2$,
where $P$ is the orthogonal projection from $\mathcal H_1$ onto $\mathcal H_1^0$.
Note that for $h\in\mathcal H_2^0$, we have $A_0R_0h=APRh=PARh$ and $R_0B_0h=PRBh$.
Clearly, $\|A_0R_0-R_0B_0\|\le\|AR-RB\|$.
Applying \rf{lip1R2} to the operators $A_0$, $B_0$, and $R_0$, we obtain
\begin{align*}
|(f(A)Ru,v)-(Rf(B)u,v)|&=|(Ru,\ov f(A)v)-(Rf(B)u,v)|\\[.2cm]
&=|(R_0u,\ov f(A_0)v)-(R_0f(B_0)u,v)|\\[.2cm]
&=|(f(A_0)R_0u,v)-(R_0f(B_0)u,v)|\\[.2cm]
&\le C\|A_0R_0-R_0B_0\|\le C\|AR-RB\|.\quad \bl
\end{align*}

\medskip

{\bf Remark.} Theorem \ref{spectra} allows us to give alternative
the proofs of \rf{ravvo}, \rf{nva} and
\rf{nva+} that do not use double operator integrals.

\begin{thm}
\label{pi2}
Let $f$ be a function defined on $\Z$.
Then
$$
\O_{f,\Z}^{\flat}(\d)=\d \|f\|_{{\rm OL}(\Z)}
$$
for $\d\in\big(0,\frac2\pi\big]$.
\end{thm}

\Pf The inequality
$$
\O_{f,\Z}^\flat(\d)\le\d\|f\|_{{\rm OL}(\Z)},\quad\d>0,
$$
is a consequence of Theorem \ref{Olip}.
Let us prove the opposite inequality for $\d\in\big(0,\frac2\pi\big]$.
Fix $\e>0$. There exists a self-adjoint operator $A$ and a bounded operator $R$ such that
$\|AR-RA\|=1$, $\s(A)\subset\Z$, and $\|f(A)R-Rf(A)\|\ge\|f\|_{{\rm OL}(\Z)}-\e$.
Put
$$
R_A\df\sum_{j\not=k}E_A(\{j\})RE_A(\{k\})=R-\sum_{j\in\Z}E_A(\{j\})RE_A(\{j\}).
$$
Clearly, $AR-RA=AR_A-R_AA$ and $f(A)R-Rf(A)=f(A)R_A-R_Af(A)$.
Thus we may assume that $R=R_A$.
Note that
$$
AR-RA=\sum_{j\not=k}(j-k)E_A(\{j\})RE_A(\{k\}).
$$
Since
$$
R=R_A=\sum_{j\not=k}\frac1{j-k}(j-k)E_A(\{j\})RE_A(\{k\}),
$$
we have $R=H\star(AR-RA)$, where
$$
H(j,k)\df\left\{\begin{array}{ll}\frac1{j-k},&\text {if}\,\,\,\,j\ne k,\\[.2cm]
0,&\text {if}\,\,\,\,j=k,
\end{array}\right.
$$
where $\star$ denotes Schur--Hadamard multiplication, see \rf{ScHad}.
It follows that
$$
\|R\|\le\|H\|_{\frak M_{\Z,\Z}}\|AR-RA\|=\|H\|_{\frak M_{\Z,\Z}}=\frac\pi2
$$
by Lemma \ref{Hmn}.

Let $\d\in\big(0,\frac2\pi\big]$. Then $\|A(\d R)-(\d R)A\|=\d$ and $\|\d R\|\le1$.
Hence,
$$
\O_{f,\Z}^\flat(\d)\ge\d\|f(A)R-Rf(A)\|\ge\d\big(\|f\|_{{\rm OL}(\Z)}-\e\big).
$$
Passing to the limit as $\e\to0$,
we obtain the desired result. $\bl$


Let $\o_{f,\frak F}$ denote the usual scalar modulus of continuity of a continuous function$f$ defined on $\frak F$.
Clearly, $\o_{f,\frak F}\le\O_{f,\frak F}$. Put $\o_f\df\o_{f,\R}$ and $\O_f\df\O_{f,\R}$.
We are going to get some estimates for the commutator modulus of continuity
$\O_{f,\frak F}^{\flat}$.
We consider first the case when $\frak F=\R$.
The following theorem is contained implicitly in \cite{NF}.

\begin{thm}
\label{kmr+}
Let $f$ be a continuous function on $\R$. Then
$$
\O_{f}^{\flat}(\d)\le2\o_f(\d/2)+2\|f(\d x)\|_{{\rm OL}(\Z)}.
$$
\end{thm}
\Pf Let $\|AR-RA\|\le\d$ with $\|R\|=1$. We can take a self-adjoint operator $A_\delta$ such that
$A_\delta A=AA_\delta$,
$\|A-A_\delta\|\le\delta/2$
and $\s(A_\delta)\subset \delta\Z$. Then
$\|f(A)-f(A_\d)\|\le\o_{f}(\d/2)$ and
\bey
\|A_\d R-RA_\d\|\le\|A_\d R-AR\|+\|AR-RA\|+\|RA-RA_\d\|\le2\d
\eey
Hence,
\begin{align*}
\|f(A)R-Rf(A)\|&\le\|f(A)R-f(A_\delta)R\|+\|f(A_\delta)R-Rf(A_\delta)\|+\|Rf(A_\d)-Rf(A)\|\\
&\le2\o_{f}(\d/2)+\|A_\delta R-RA_\delta\|\cdot\|f\|_{{\rm OL}(\d\Z)}
\le2\o_{f}(\d/2)+2\delta\|f\|_{{\rm OL}(\d\Z)}\\
&=2\o_{f}(\d/2)+2\|f(\d x)\|_{{\rm OL}(\Z)}.\quad\bl
\end{align*}
\begin{thm}
\label{kmr-}
Let $f$ be a continuous function on $\R$. Then
$$
\O_f^{\flat}(\d)\ge
\max\left\{\o_f(\d),\frac2\pi\|f(\d x)\|_{{\rm OL}(\Z)}\right\}
$$
for all $\d>0$.
\end{thm}

\Pf Clearly, $\o_f\le\O_f\le\O_f^{\flat}$. It remains to prove that
$\|f(\d x)\|_{{\rm OL}(\Z)}\le\frac\pi2\O_f^{\flat}(\d)$.
We have
$$
\O_f^{\flat}(\d)\ge\O_{f,\d\Z}^{\flat}(\d)=\O_{f(\d x),\Z}^{\flat}(1)\ge\O_{f(\d x),\Z}^{\flat}\Big(\frac2\pi\Big)
=\frac2\pi\big\|f(\d x)\big\|_{{\rm OL}(\Z)}
$$
by Theorem \ref{pi2}. $\bl$

We consider now similar estimates of $\O_{f,\frak F}^{\flat}$ for an arbitrary closed subset $\fF$ of $\R$.
Recall that a subset $\L$ of $\R$ is called a $\d$ net for $\fF$ if
$\frak F\subset\bigcup\limits_{t\in\L}[t-\d,t+\d]$.

\begin{thm}
\label{kme+}
Let $f$ be a continuous function on a closed subset $\frak F$ of $\R$. Suppose that $\fF_\d$ is a subset of $\fF$ that forms a $\d/2$ net of $\fF$. Then
$$
\O_{f,\frak F}^{\flat}(\d)\le2\o_{f,\frak F}(\d/2)+2\d\|f\|_{{\rm OL}(\frak F_\d)}.
$$
\end{thm}

\Pf The proof is similar to the proof of Theorem \ref{kmr+}.
It suffices to replace the $\d/2$ net $\d\Z$ of $\R$ with
the $\d/2$ net $\fF_\d$ of $\fF$. $\bl$

\begin{thm}
\label{kme-}
Let $f$ be a continuous function on a closed subset $\frak F$ of $\R$ and let $\d>0$.
Suppose that $\L$ and $\rm M$ are closed subsets of $\frak F$ such that
$(\L-{\rm M})\cap(-\d,\d)\subset\{0\}$.
Then
$$
\O_{f,\frak F}^{\flat}(\d)\ge
\max\left\{\o_{f,\fF}(\d),\frac\d2\|\dg_0 f\|_{\fM_{\L,\rm M}}\right\}.
$$
\end{thm}

\Pf Clearly, $\o_{f,\frak F}\le\O_{f,\frak F}\le\O_{f,\frak F}^{\flat}$.
Note that
$$
\|\dg_0 f\|_{\fM_{\L,\rm M}}=\sup_{a>0}\|\dg_0 f\|_{\fM_{\L\cap[-a,a]),{\rm M}\cap[-a,a]}}.
$$
Thus it suffices to prove that
$$
\O_{f,\frak F}^{\flat}(\d)\ge\frac\d2\|\dg_0 f\|_{\frak M_{\L,\rm M}}
$$
in the case when $\L$ and $\rm M$ are bounded.

Let $\e>0$. There exist positive regular Borel measures $\l$ on $\L$, $\mu$
on $\rm M$, and a function $k$ in $L^2(\L\times{\rm M},\l\otimes\mu)$ such that
$\|k\|_{\mB_{\L,{\rm M}}^{\l,\mu}}=1$
and $\|k\dg_0 f\|_{\mB_{\L,{\rm M}}^{\l,\mu}}\ge\|\dg_0 f\|_{\fM_{\L,{\rm M}}}-\e$.
We define the function $k_0$ in $L^2(\L\times{\rm M},\l\otimes\mu)$ by
$$
k_0(x,y)\df\left\{\begin{array}{ll}k(x,y),&\text {if}\,\,\,\,x\ne y,\\[.2cm]
0,&\text {if}\,\,\,\,x=y.
\end{array}\right.
$$
Then $k\dg_0 f=k_0\dg_0 f$ and
$\|k_0\|_{\mB_{\L,{\rm M}}^{\l,\mu}}\le2$.
Put $\Phi(x,y)\df f_\d(x-y)$ where $f_\d$ denotes the same as in Corollary \ref{Hxy}.
We define the self-adjoint operators $A:L^2(\L,\l)\to L^2(\L,\l)$ and
$B:L^2({\rm M},\mu)\to L^2({\rm M},\mu)$ by $(Af)(x)\df xf(x)$ and
$(Bg)(y)\df yg(y)$.
Put
$$
h(x,y)\df\Phi(x,y)k(x,y)=\Phi(x,y)k_0(x,y).
$$
Clearly,
$$
\|h\|_{\mB_{\L,{\rm M}}^{\l,\mu}}
\le\|\Phi\|_{\fM_{\L,{\rm M}}^{\l,\mu}}\|k\|_{\mB_{\L,{\rm M}}^{\l,\mu}}
\le\|\Phi\|_{\fM_{\R,\R}}\le\frac2{\d}
$$
by Corollary \ref{Hxy}.

Clearly, $A\mI_h-\mI_h B=\mI_{k_0}$ and
$f(A)\mI_h-\mI_h f(B)=\mI_{k_0\dg_0 f}$.
(Recall that $\mI_\f$ is the integral operator from $L^2({\rm M},\mu)$ into $L^2(\L,\l)$ with kernel $\f\in L^2(\L\times{\rm M},\l\otimes\nu)$.)
Then
$$
\left\|\frac\d2\mI_h\right\|=\frac\d2\|h\|_{\mB_{\L,{\rm M}}^{\l,\mu}}\le1,
$$
$$
\left\|A\left(\frac\d2\mI_h\right)-\left(\frac\d2\mI_h\right) B\right\|=\frac\d2\|k_0\|
_{\mB_{\L,{\rm M}}^{\l,\mu}}\le\d,
$$
and
$$
\left\|f(A)\left(\frac\d2\mI_h\right)-\left(\frac\d2\mI_h\right) f(B)\right\|
=\frac\d2\|k_0\dg_0 f\|_{\mB_{\L,{\rm M}}^{\l,\mu}}
\ge
\frac\d2\big(\|\dg_0 f\|_{\fM_{\L,{\rm M}}}-\e\big).
$$
Hence, $\O_{f,\frak F}^{\flat}(\d)\ge\frac\d2\big(\|\dg_0 f\|_{\fM_{\L,{\rm M}}}-\e\big)$
for every $\e>0$. $\bl$

Theorem \ref{kme-} allows us to obtain another proof of Theorem 4.17 in \cite{AP4}.

\begin{thm}
\label{JW}
Let $f$ be a continuous function on an unbounded closed subset $\fF$ of $\R$.
Suppose that $\O_{f, \fF}(\d)<\be$ for $\d>0$.
Then
the function $t\mapsto t^{-1}f(t)$ has a finite limit as $|t|\to\be$, $t\in \fF$.
\end{thm}

\Pf
Assume the contrary. Then there exists a sequence $\{\l_n\}_{n=1}^\be$ in $\fF$ such that
$|\l_{n+1}|-|\l_n|>1$ for all $n\ge1$, $\lim_{n\to\be}|\l_n|=\be$ and the sequence $\{\l_n^{-1}f(\l_n)\}_{n=1}^\be$
has no finite limit. Denote by $\L$ the image of the sequence $\{\l_n\}_{n=1}^\be$.
Then $\|f\|_{{{\rm OL}(\L)}}=\be$. This fact is contained implicitly in \cite{JW}.
Indeed, Theorem 4.1 in \cite{JW} implies that every operator Lipschitz function $f$ is differentiable
at every non-isolated point. It is well known that the same argument gives us the differentiability
at $\be$ in the following sense: the function $t\mapsto t^{-1}f(t)$ has a finite limit as $|t|\to\be$,
provided the domain of $f$ is unbounded. Applying Theorem \ref{kme-}
for $M=\L$ and $\d=1$, we find that $\O_{f, \fF}(1)=\be$. $\bl$

We need the following known result, see \cite{KST}.
We give the proof for the reader's convenience.

\begin{thm}
\label{kusL}
Let $f$ be a bounded continuous function on a closed subset
$\frak F$ of $\R$. Suppose that $f\in{\rm OL}\big((-\be,1]\cap\fF\big)$
and $f\in{\rm OL}\big([-1,\be\big)\cap\frak F)$. Then $f\in{\rm OL}(\frak F)$
and
$$
\|f\|_{{\rm OL}(\frak F)}\le
C\big(\|f\|_{{\rm OL}((-\be,1]\cap\fF)}+\|f\|_{{\rm OL}([-1,\be)\cap\fF)}
+\sup_\fF|f|\big),
$$
where $C$ is a numerical constant.
\end{thm}

\Pf Put $\frak F_1\df\fF\cap(-\be,-1]$, $\frak F_2\df\frak F\cap[-1,1]$,
and $\frak F_3\df \fF\cap[1,\be)$.
We have
\begin{align*}
\|f\|_{{\rm OL}(\frak F)}\le\|\dg_0f\|_{{\fM}_{\frak F,\frak F}}
\le&\sum_{j=1}^3\sum_{k=1}^3\|\dg_0f\|_{{\fM}_{\frak F_j,\frak F_k}}
=\sum_{j=1}^3\|\dg_0f\|_{{\fM}_{\frak F_j,\frak F_j}}\\[.2cm]
&+
2\|\dg_0f\|_{{\fM}_{\frak F_1,\frak F_2}}+2\|\dg_0f\|_{{\fM}_{\frak F_2,\frak F_3}}
+2\|\dg_0f\|_{{\fM}_{\frak F_1,\frak F_3}}.
\end{align*}

Each term $\|\dg_0f\|_{{\fM}_{\frak F_j,\frak F_k}}$ except $\|\dg_0f\|_{{\fM}_{\frak F_1,\frak F_3}}$
can be estimated in terms of \lb$2\|f\|_{{\rm OL}(\frak F_1\cup\frak F_2)}$ or
 $2\|f\|_{{\rm OL}(\frak F_2\cup\frak F_3)}$.

Let us estimate $\|\dg_0f\|_{{\fM}_{\frak F_1,\frak F_3}}$.
We have
\begin{align*}
\|\dg_0f\|_{{\fM}_{\frak F_1,\frak F_3}}&=\left\|\frac{f(x)-f(y)}{x-y}\right\|_{{\fM}_{\frak F_1,\frak F_3}}\\[.2cm]
&\le\left(\sup_{\frak F_1}|f|\right)\left\|\frac{1}{x-y}\right\|_{{\fM}_{\frak F_1,\frak F_3}}+
\left(\sup_{\frak F_3}|f|\right)\left\|\frac{1}{x-y}\right\|_{{\fM}_{\frak F_1,\frak F_3}}\\[.2cm]
&\le2\left(\sup_\frak F|f|\right)\left\|\frac{1}{x-y}\right\|_{{\fM}_{\frak F_1,\frak F_3}}
\le2\sup_\frak F|f|
\end{align*}
because by Corollary \ref{Hxy},
$$
\left\|\frac{1}{x-y}\right\|_{{\fM}_{\frak F_1,\frak F_3}}\le\|f_2(x-y)\|_{{\fM}_{\R,\R}}\le1,
$$
where $f_2$ means the same as in Corollary \ref{Hxy}.

Thus
$$
\|f\|_{{\rm OL}(\frak F)}\le6\|f\|_{{\rm OL}(\frak F_1\cup\frak F_2)}+
4\|f\|_{{\rm OL}(\frak F_2\cup\frak F_3)}+4\sup_\frak F|f|.\quad\bl
$$

\

\section{\bf The operator Lipschitz norm of the function $\bs{x\mapsto|x|}$
on subsets of $\bs{\R}$}
\setcounter{equation}{0}
\label{|x|}

\

In this section we obtain sharp estimates of the operator modulus of continuity of the function $x\mapsto|x|$ on certain subsets of the real line. This allows us to obtain sharp estimates of $\big\|\,|S|-|T|\,\big\|$ for arbitrary bounded linear operators $S$ and $T$. Note that our estimates considerably improve earlier results of \cite{Ka}.

Put $\Abs(x)\df|x|$. For $J\subset[0,\be)$, we put
$\log (J)\df\{\log t:t\in J,\,t>0\}$.

\begin{thm}
\label{lipm}
There exist positive numbers $C_1$ and $C_2$ such that
\bey
C_1\log\big(2+|\log(J_1\cap J_2)|\big)\le\|\Abs\|_{{\rm OL}((-J_1)\cup J_2)}
\le C_2\log\big(2+|\log(J_1\cap J_2)|\big)
\eey
for all intervals $J_1$ and $J_2$ in $(0,\be)$.
\end{thm}

\Pf Put $J=J_1\cap J_2$. Let us first establish the lower estimate.
Note that $\|\Abs\|_{{\rm OL}((-J_1)\cup J_2)}\ge\|\Abs\|_{{\rm OL}(J_2)}=1$.
This proves the lower estimate in the case $|\log(J)|\le1$.
In the case $|\log(J)|>1$ we have
\begin{align*}
\|\Abs\|_{{\rm OL}((-J_1)\cup J_2)}&\ge\|\Abs\|_{{\rm OL}((-J)\cup J)}\ge
\left\|\frac{|x|-|y|}{x-y}\right\|_{\fM_{-J,J}}=\left\|\frac{x-y}{x+y}\right\|_{\frak M_{J,J}}\\[.2cm]
&\ge c_1\log(1+|\log(J)|)\ge c_2\log(2+|\log(J)|)
\end{align*}
by Theorem \ref{albe-}.

We proceed now to the upper estimate.
We consider first the case when $J=J_1$. Then
$$
\|\Abs\|_{{\rm OL}((-J_1)\cup J_2)}\le\|\Abs\|_{{\rm OL}((-J_1)\cup [0,\be))}\le
2+2\left\|\frac{x-y}{x+y}\right\|_{\frak M_{J,[0,\be)}}
$$
and we can apply Theorem \ref{albe+}. The case $J=J_2$ is similar.
Suppose that $J\ne J_1$ and $J\ne J_2$. Then
$\inf J_1\ne\inf J_2$. Let $\inf J_1>\inf J_2$.
Put $a\df\inf J_1$ and $b\df\sup J_2$. Then
$$
\|\Abs\|_{{\rm OL}((-J_1)\cup J_2)}\le\|\Abs\|_{{\rm OL}((-\be_-a]\cup [0,b))}\le
2+2\left\|\frac{x-y}{x+y}\right\|_{\frak M_{[a,\be),[0,b)}}
$$
and the result follows from Theorem \ref{albe}. $\bl$

Let us state two special cases of Theorem \ref{lipm}.

\begin{thm}
\label{lipm1}
There exist positive constants $C_1$ and $C_2$ such that
$$
C_1\log(2+\log(ba^{-1}))\le\|\Abs\|_{{\rm OL}((-\be,0]\cup[a,b])}
\le C_2\log(2+\log(ba^{-1}))
$$
for every $a,b\in(0,\be)$ with $a<b$.
\end{thm}

\begin{thm}
\label{lipm2}
There exist positive constants $C_1$ and $C_2$ such that
$$
C_1\log(2+\log_+(ba^{-1}))\le\|\Abs\|_{{\rm OL}((-b,0]\cup[a,\be))}
\le C_2\log(2+\log_+(ba^{-1}))
$$
for every $a,b\in(0,\be)$.
\end{thm}

\begin{thm}
\label{fa}
Let $\xi_a=\Abs\big|{[-a,\be)}$ and $\eta_a=\Abs\big|{[-a,a]}$, where $a>0$.
Then there exist positive numbers $C_1$ and $C_2$ such that
$$
C_1\d\log(2+\log_+(a\d^{-1}))\le\O_{\eta_a}(\d)\le\O_{\xi_a}(\d)\le C_2\d\log(2+\log_+(a\d^{-1}))
$$
for every $\d>0$.
\end{thm}

\Pf Put $\frak F_\d\df[-a,\be)\setminus(0,\d)$. Clearly, $\frak F_\d$
is a $\d$-net of $(-\be,a]$. Hence, by Theorem \ref{kme+} we have
$$
\O_{\xi_a}(\d)\le\O_{\xi_a}^{\flat}(\d)\le\d+2\d\|\xi_a\|_{{\rm OL}(\frak F_\d)}.
$$
Applying Theorem \ref{lipm2}, we obtain the desired upper estimate.

To obtain the lower estimate, we use Theorem \ref{kme-}. Clearly,
$\O_{\eta_a}(\d)\ge\d$ for all $\d\in(0,a]$. Thus it suffices to consider the case
$\d\in(0,\frac a2)$. Put
$\L=[-a,0]$ and ${\rm M}=[\d,a]$. By Theorem \ref{kme-},
$$
\O_{\eta_{a}}(\d)\ge\frac12\O_{\eta_{a}}^{\flat}(\d)\ge\frac\d4\|\dg_0\eta_a\|_{\frak M_{\L,\rm M}}.
$$
It remains to apply Theorem \ref{albe-}. $\bl$

\begin{thm}
\label{sverkhu}
There exists a positive number $C$ such that
$$
\big\|\,|A|-|B|\,\big\|\le C\|A-B\|\log\left(2+\log\frac{\|A\|+\|B\|}{\|A-B\|}\right)
$$
for every bounded self-adjoint operators $A$ and $B$.
\end{thm}

\Pf This is a special case of Theorem \ref{fa} that corresponds to $a=\|A\|+\|B\|$. $\bl$

Theorem also \ref{fa} allows us to prove that the upper estimate in Theorem \ref{sverkhu} is sharp.

\begin{thm}
\label{snizu}
Let $a>0$. There is a positive number $c$ such that for every $\d\in(0,a)$, there exist self-adjoint operators $A$ and $B$ such that $\|A+B\|\le a$, $\|A-B\|\le\d$, but
$$
\big\|\,|A|-|B|\,\big\|\ge\,c\d\log\left(2+\log\frac{a}{\d}\right).
$$
\end{thm}

We proceed now to the case of arbitrary (not necessarily self-adjoint) operators. Recall that
for a bounded operator $S$ on Hilbert space, its modulus $S$ is defined by
$$
|S|\df(S^*S)^{1/2}.
$$

\begin{thm}
\label{Kato}
There exists a positive number $C$ such that
$$
\big\|\,|S|-|T|\,\big\|\le C\|S-T\|\log\left(2+\log\frac{\|S\|+\|T\|}{\|S-T\|}\right)
$$
for every bounded operators $S$ and $T$.
\end{thm}

\Pf
Put
$$
A=\left(\begin{matrix}\0& S^*\\[.2cm]S&\0\end{matrix}\right)\quad\mbox{and}\quad
B=\left(\begin{matrix}\0&T^*\\[.2cm]T&\0\end{matrix}\right).
$$
Clearly, $A$ and $B$ are self-adjoint operators with
$$
|A|=\left(\begin{matrix}|S|&\0\\[.2cm]\0&|S^*|\end{matrix}\right)\quad\mbox{and}\quad
|B|=\left(\begin{matrix}|T|&\0\\[.2cm]\0&|T^*|\end{matrix}\right).
$$
Hence,
\begin{align*}
\big\|\,|S|-|T|\,\big\|&\le\big\|\,|A|-|B|\,\big\|\le C\|A-B\|\log\left(2+\log\frac{\|A\|+\|B\|}{\|A-B\|}\right)\\[.2cm]
&=C\|S-T\|\log\left(2+\log\frac{\|S\|+\|T\|}{\|S-T\|}\right).\quad \bl
\end{align*}

\medskip

{\bf Remark.} Theorem \ref{Kato} significantly improves Kato's inequality obtained in \cite{Ka}:
$$
\big\|\,|S|-|T|\,\big\|\le \frac1\pi\|S-T\|\left(2+\log\frac{\|S\|+\|T\|}{\|S-T\|}\right).
$$

\

\section{\bf The operator modulus of continuity of a certain piecewise linear function}
\setcounter{equation}{0}
\label{kuslin}

\

In this section we obtain a sharp estimate for the operator modulus of continuity of
the piecewise linear function $\vk$ defined by
$$
\vk(t)\df\left\{\begin{array}{ll}1,&\text {if}\,\,\,\,t\ge1,\\[.2cm]
t,&\text {if}\,\,\,-1<t\le1.\\[.2cm]
-1,&\text {if}\,\,\,\,t>1.
\end{array}\right.
$$
The results obtained in this section will be used in the next section to estimate the operator modulus of continuity of functions concave on $\R_+$.

It is easy to see that $\vk(t)=\frac12\big(|1+t|-|1-t|\big)$.

\begin{thm}
\label{lipm4}
There exist positive numbers $C_1$ and $C_2$ such that
$$
C_1\log|\log\d|\le\|\vk\|_{{\rm OL}((-\be,-1-\d]\cup[-1,1]\cup[1+\d,\be))}\le C_2\log|\log\d|
$$
for every $\d\in(0,\frac12)$.
\end{thm}

\Pf Put $\vk_1=\vk\big|{((-\be,-1-\d]\cup[-1,1])}$ and
$\vk_2=\vk\big|{([-1,1]\cup[1+\d,\be))}$. Note that
$$
\vk_1(t)=\frac12\big(|1+t|-1+t\big)
\quad\mbox{and}\quad
\vk_2(t)=\frac12\big(1+t-|t-1|\big).
$$
It follows from Theorem \ref{lipm2} that
$$
C_1\log|\log\d|\le\|\vk_1\|_{{\rm OL}}\le C_2\log|\log\d|
$$
and
$$
C_1\log|\log\d|\le\|\vk_2\|_{{\rm OL}}\le C_2\log|\log\d|.
$$
Thus the desired lower estimate is evident and
the required upper estimate follows from Theorem \ref{kusL}. $\bl$

\begin{thm}
\label{skok}
There exist positive numbers $c_1$ and $c_2$ such that
$$
c_1\d\log\big(1+\log(1+\d^{-1})\big)\le\O_\vk(\d)\le c_2\d\log\big(1+\log(1+\d^{-1})\big)
$$
for every $\d>0$.
\end{thm}

\Pf  Note that $\lim\limits_{t\to\be}t\log\big(1+\log(1+t^{-1})\big)=1$.
Thus it suffices to consider the case when $0<\d\le\frac12$.
Put $\frak F_\d\df(-\be,-1-\d]\cup[-1,1]\cup[1+\d,\be)$. Clearly, $\frak F_\d$
is a $\d$-net for $\R$. Hence, by Theorem \ref{kme+}, we have
$$
\O_{\vk}(\d)\le\O_{\vk}^{\flat}(\d)\le\d+2\d\|\vk\|_{{\rm OL}(\frak F_\d)}.
$$
The desired upper estimate follows now from Theorem \ref{lipm4}.

To obtain the lower estimate we can apply Theorem \ref{fa}
because $\vk(t)=\frac12(|1+t|-1+t)$ for $t\le1$. $\bl$

\

\section{\bf Operator moduli of continuity of concave functions on $\bs{\R_+}$.}
\setcounter{equation}{0}
\label{conc}
\

Recall that in \cite{AP2} we proved that if $f$ is a continuous function on $\R$, then its operator modulus of continuity $\O_f$ admits the estimate
$$
\O_f(\d)\le\const\,\d\int_\d^\be\frac{\o_f(t)}{t^2}\,dt=\const\int_1^\infty\frac{\o_f(t\d)}{t^{2}}ds,\quad\d>0.
$$

In this section we show that if $f$ vanishes on $(-\be,0]$ and is a concave nondecreasing function on $[0,\be)$, then the above estimate can be considerably improved.

We also obtain several other estimates of operator moduli of continuity.

\begin{thm}
\label{fM}
Suppose that $f^{\prime\prime}=\mu\in\M(\R)$ (in the distributional sense),
$\mu(\R)=0$, and
$$
\int_\R\log(\log(|t|+3))\,d|\mu|(t)<\be.
$$
Then
$$
\O_f(\d)\le c\,\|\mu\|_{\M(\R)}\d\log(\log(\d^{-1}+3)),
$$
where $c$ is a numerical constant.
\end{thm}

\Pf Put
\bay
\label{fis}
\f_s(t)\df\frac12\big(|t|+|s|\big)-\frac{|t-s|}2,\quad s,t\in\R.
\ey
It is easy to see that
$$
\f_s(t)\df\frac{|s|}2\vk\left(\frac{2t}s-1\right)+\frac{|s|}2\quad\mbox{for}\quad s\ne0.
$$
Clearly,
\bay
\label{vtoraya}
\f_s^{\prime\prime}=\d_0-\d_s\quad\mbox{and}\quad\f_s(0)=0.
\ey
Theorem \ref{skok} implies that
\begin{align}
\label{modcon}
\O_{\f_s}(t)&\le\const\,t\,\log\left(1+\log\left(1+\frac{|s|}{2t}\right)\right)\nonumber\\[.2cm]
&\le\const\,t\,\log\left(1+\log\left(1+\frac{|s|}{t}\right)\right),\quad t>0.
\end{align}
It is easy to see that
$$
t\,\log\Big(1+\log\big(1+t^{-1}|s|\big)\Big)\le\const\big(\log(\log(|s|+3))\big)\,t\,
\log\Big(\log\big(t^{-1}+3\big)\Big).
$$
To complete the proof, it suffices to observe that
$$
f(t)=at +b-\int_\R\f_s(t)\,d\mu(s),\quad\mbox{for some}\quad a,~b\in\C,
$$
which follows easily from \rf{vtoraya}. $\bl$

The assumption that $\mu(\R)=0$ in the hypotheses of Theorem \ref{fM} is essential.
Moreover, the following result holds.

\begin{thm}
\label{sushch}
Suppose that $f^{\prime\prime}=\mu\in\M(\R)$ and
$\mu(\R)\not=0$. Then $\O_f(t)=\be$ for every $t>0$.
\end{thm}

\Pf
Indeed, it is easy to see that
$f^\prime(t)=c+\mu((-\be,t))$ for almost all $t\in\R$.
Hence,
$$
\lim_{t\to\be}\frac{f(t)}t=\lim_{t\to\be}f^{\prime}(t)=c+\mu(\R)\quad\text{and}\quad
\lim_{t\to-\be}\frac{f(t)}t=\lim_{t\to-\be}f^{\prime}(t)=c.
$$
The result follows from Theorem \ref{JW}. $\bl$

Let $G$ be an open subset of $\R$. Denote by $\M_{\rm loc}(G)$ the set of all distributions
on $G$ that are locally (complex) measures.

\begin{thm}
\label{log}
Let $f\in C(\R)$. Put $\mu\df f^{\prime\prime}$ in the sense of distributions.
Suppose that $\lim\limits_{|t|\to\be}t^{-1}f(t)=0$, $\mu\big|(\R\setminus\{0\})\in\M_{\rm loc}(\R\setminus\{0\})$
and
$$
\int_{\R\setminus\{0\}}\log(1+\log(1+|s|))\,d|\mu|(s)<\be.
$$
Then
$$
\O_f(\d)\le\const\,\d\int_{\R\setminus\{0\}}\log\Big(1+\log\big(1+|s|\d^{-1}\big)\Big)\,d|\mu|(s).
$$
\end{thm}

\Pf Put
$$
g(t)=-\int_{\R\setminus\{0\}}\f_s(t)\,d\mu(s),
$$
where $\f_s$ is defined by \rf{fis}.
Inequality \rf{modcon} implies that
\bay
\label{Om}
\O_g(\d)\le\const\,\d\int_{\R\setminus\{0\}}\log\Big(1+\log\big(1+|s|\d^{-1}\big)\Big)\,d|\mu|(s).
\ey
In particular, $g$ is continuous on $\R$. Clearly, $g^{\prime\prime}=f^{\prime\prime}$
on $\R\setminus\{0\}$. Hence, $f(x)-g(x)=a|x|+bx+c$ for some $a,\,b,\,c\in\C$. It follows from \rf{Om} that
$$
\lim_{|t|\to\be}\left|\frac{g(t)}t\right|\le\lim_{t\to\be}\frac{\o_g(t)}t\le\lim_{t\to\be}\frac{\O_g(t)}t=0=\lim_{|t|\to\be}\frac{f(t)}t
$$
which implies that $f-g=\const$. $\bl$

\begin{cor}
Let $a>0$ and let $f$ be a continuous function on $\R$ that is constant on $\R\setminus(-a,a)$.
Put $\mu\df f^{\prime\prime}$ in the sense of distributions.
Suppose that $\mu\big|(\R\setminus\{0\})\in\M_{\rm loc}(\R\setminus\{0\})$ and
\bay
\label{s2s}
C\df\sup_{s>0}|\mu|\big([s,2s]\cup[-2s,-s]\big)<\be.
\ey
Then
$$
\O_f(\d)\le C\const\,\d\,\Big(\log \frac a\d\Big)\log\Big(\log \frac a\d\Big)\,\quad\mbox{for}\quad
\d\in\Big(0,\frac a3\Big).
$$
\end{cor}

\Pf By Theorem \ref{log},
\begin{align*}
\O_f(\d)\le&\const\,\d\!\left(\int_{0}^a\!\log\big(1\!+\!\log\big(1\!+\!s\d^{-1}\big)\big)\,d|\mu(s)|\!+
\!\int_{0}^a\!\log\big(1\!+\!\log\big(1\!+\!s\d^{-1}\big)\big)\,d|\mu(-s)|\!\right)\\[.2cm]
=&\const\,\d\sum_{n\ge0}\int_{2^{-n-1}a}^{2^{-n}a}\log\big(1+\log\big(1+s\d^{-1}\big)\big)\, d|\mu|(s)\\[.2cm]
&+
\const\,\d\sum_{n\ge0}
\int_{2^{-n-1}a}^{2^{-n}a}\log\big(1+\log\big(1+s\d^{-1}\big)\big)\,d|\mu|(-s).\end{align*}
It follows now from \rf{s2s} and the inequality
$$
\log(1+\log(1+\a x))\le2\log(1+\log(1+x)), \quad 0<x<\be,\quad 1<\a\le2,
$$
that
\begin{align*}
\O_f(\d)\le&\const\,\d\sum_{n\ge0}\int_{2^{-n-1}a}^{2^{-n}a}
\log\big(1+\log\big(1+s\d^{-1}\big)\big)\frac{ds}s\\[.2cm]
=&\const\,\d\int_{0}^a\log\big(1+\log\big(1+s\d^{-1}\big)\big)\frac{ds}s\\[.2cm]
=&\const\,\d\int_{0}^{a/\d}\log\big(1+\log(1+s)\big)\frac{ds}s\\[.2cm]
\le&\const\,\d+ \const\,\d\int_{1}^{a/\d}\log\big(1+\log(1+s)\big)\frac{ds}s\\[.2cm]
=& \const\,\d\left(1+\big(\!\log\big(1+\log(1+s)\big)\log s\big)\Big|_1^{a/\d}
\!-\int_{1}^{a/\d}\!\frac{\log s\,ds}{(1+s)\log\big(1+\log(1+s)\big)}\right)\\[.2cm]
\le& \const\,\d+ \const\,\d\big(\log(1+\log(1+s)\big)\log s\big)\Big|_1^{a/\d}
\le \const\,\d\Big(\log \frac a\d\Big)\log\Big(\log \frac a\d\Big)
\end{align*}
for sufficiently small $\d$. $\bl$

\begin{cor}
Let $f$ be a continuous function on $\R$ that is constant on $\R\setminus(-a,a)$. Suppose that $f$ is twice differentiable on $\R\setminus\{0\}$    and
$$
C\df\sup_{s\ne0}\big|sf^{\prime\prime}(s)\big|<\be.
$$
Then
$$
\O_f(\d)\le \const\,C\,\d\,\Big(\log \frac a\d\Big)\log\Big(\log \frac a\d\Big)\,\quad\mbox{for}\quad
\d\in\Big(0,\frac a3\Big).
$$
\end{cor}

The following result shows that in a sense Theorem \ref{fM} cannot be improved.

\begin{thm}
\label{sh}
Let $h$ be a positive continuous function on $\R$.
Suppose that for every $f\in C(\R)$ such that
$$
f^{\prime\prime}=\mu\in\M(\R),\quad
\mu(\R)=0,\quad\mbox{and}\quad \int_\R h(t)\,d|\mu|(t)<\be,
$$
we have $\O_f(\d)<\be$, $\d>0$. Then for some positive number $c$,
$$
h(t)\ge c\log(\log(|t|+3)),\quad t\in\R.
$$
\end{thm}

We need the following lemma, in which $\f_s$ is the function defined by \rf{fis}.

\begin{lem}
\label{s10}
There is a positive number $c$ such that for every
$s\ge10$, there exist self-adjoint operators
$A$ and $B$ satisfying the conditions:
$$
\s(A),~\s(B)\subset\left(\frac s2,\frac{3s}2\right),\quad
\|A-B\|\le1,\quad\mbox{and}\quad \|\f_s(A)-\f_s(B)\|\ge c\log\log s.
$$
\end{lem}

\Pf Clearly, it suffices to prove the lemma for sufficiently large $s$. By Theorem \ref{fa},
there exist self-adjoint operators
$A_0$ and $B_0$ such that $\|A_0\|, \|B_0\|<1$, $\|A_0-B_0\|\le2/s$, and
$\big\|\,|A_0|-|B_0|\,\big\|\ge\const s^{-1}\log\big(2+\log s\big)$.
Put $A\df sI+\frac s2 A_0$ and $B\df sI+\frac s2 B_0$.
Then $\s(A),\,\s(B)\subset\big(\frac s2,\frac{3s}2\big)$ and
$\|A-B\|\le1$. Let us estimate $\|\f_s(A)-\f_s(B)\|$.
Clearly,
$$
\f_s(A)\!-\!\f_s(B)=\frac s4(A_0-B_0)-\frac s4(|A_0|-|B_0|).
$$
Hence,
\begin{align*}
\|\f_s(A)-\f_s(B)\|&\ge\frac s4\big\|\,|A_0|-|B_0|\,\big\|-\frac s4\|A_0-B_0\|\\[.2cm]
&\ge\const\log\log s-\frac12\ge\const\log\log s
\end{align*}
for sufficiently large $s$. $\bl$

\medskip

{\bf Proof of Theorem \ref{sh}.}
Assume the contrary. Then there exists a sequence $\{s_n\}$ of real numbers
such that $\lim\limits_{n\to\be}|s_n|=\be$ and $\lim\limits_{n\to\be}(\log(\log(|s_n|)))^{-1}h(s_n)=0$. Passing to a subsequence, we can reduce the situation to the case when
$s_n>0$ for all $n$ or $s_n<0$ for all $n$. Without loss of generality we may assume that $s_n>0$
for all $n$. Moreover, we may also assume that $s_1\ge10$, $s_{n+1}\ge2 s_n$
and $\log\log s_n\ge n^3(1+h(s_n))$ for every $n\ge1$. Put $\a_n\df n(\log\log s_n)^{-1}$ for $n\ge1$ and
$f(t)\df\sum\limits_{n\ge1}\a_n\f_{s_n}(t)$. Note that the series converges for every $t$ because
$\s\df\sum\limits_{n\ge1}\a_n<\be$. Moreover,
$$
f^{\prime\prime}=\s\d_0-\sum\limits_{n\ge1}\a_n\d_{s_n}\quad\mbox{and}\quad
\s h(0)+\sum\limits_{n\ge1}\a_n h(s_n)<\be.
$$
By Lemma \ref{s10}, there exist two sequences $\{A_n\}_{n\ge1}$ and $\{B_n\}_{n\ge1}$
of self-adjoint operators such that
$$
\s(A_n),~\s(B_n)\subset\left(\frac{s_n}2,\frac{3s_n}2\right),\quad
\|A_n-B_n\|\le1,
$$
and
$$
\|\f_{s_n}(A_n)-\f_{s_n}(B_n)\|\ge c\log\log s_n.
$$
Note that $\f_{s_k}(A_n)=\f_{s_k}(B_n)=s_kI$ for $k<n$. Also, $\f_{s_k}(A_n)=A_n$
and $\f_{s_k}(B_n)=B_n$ for $k>n$.
Hence,
$$
f(A_n)-f(B_n)=\a_n(\f_{s_n}(A_n)-\f_{s_n}(B_n))+\sum_{k>n}\a_k(A_n-B_n),
$$
and so
\bey
\|f(A_n)-f(B_n)\|\ge\a_n\big\|\f_{s_n}(A_n)-\f_{s_n}(B_n)\big\|-\sum_{k>n}\a_k\|A_n-B_n\|\\[.2cm]
\ge
C\a_n\log\log s_n-\sum_{k>n}\a_k\to \be\quad \text {as} \quad n\to\be.
\eey
Thus $\O_f(1)=\be$ and we get a contradiction. $\bl$


In \cite{AP2} it was proved that
$$
\O_f(\d)\le\int_1^\be\frac{\o_f(\d s)}{s^2}\,ds
$$
for every $f\in C(\R)$. The following theorem shows that this estimate
can be improved essentially for functions $f$ concave on a ray.

\begin{thm}
\label{ostar}
Let $f$ be a continuous nondecreasing function such that
$f(t)=0$ for $t\le0$, $\lim\limits_{t\to\be}t^{-1}f(t)=0$, and $f$ is concave on $[0,\be)$. Then
$$
\O_f(\d)\le c\int_e^\be\frac{f(\d s)\,ds}{s^2\log s},
$$
where $c$ as a numerical constant.
\end{thm}

\Pf Let $\mu=-f^{\prime\prime}$ (in the distributional sense). Clearly, $\mu=0$ on $(-\be,0)$ and
$\mu$ is a positive regular measure on $(0,\be)$ because $f$ is concave on $(0,\be)$.
Hence, $\mu\in\M_{\rm{loc}}(\R\setminus\{0\})$.
By Theorem \ref{log}, we have
$$
\O_f(\d)
\le\const\,\d\int_0^\be\log(1+\log(1+s\d^{-1}))\,d\mu(s).
$$
To estimate this integral, we use the equality $f^{\prime}(t)=\mu(t,\be)$ for almost all $t>0$
and apply the Tonelli theorem twice.
\begin{align*}
\d\int_0^\be\log(1+\log(1+&s\d^{-1}))\,d\mu(s)
=\int_0^\be\left(\int_0^s\frac{dt}{(1+\log(1+t\d^{-1}))(1+t\d^{-1})}\right)\,d\mu(s)\\[.2cm]
=&\int_0^\be\frac{f^\prime(t)\,dt}{(1+\log(1+t\d^{-1}))(1+t\d^{-1})}\\[.2cm]
=&\,\d^{-1}\int_0^\be\left(\int_t^\be\frac{(2+\log(1+s\d^{-1}))\,ds}{(1+
\log(1+s\d^{-1}))^2(1+s\d^{-1})^2}\right)f^\prime(t)\,dt\\[.2cm]
=&\,\d^{-1}\int_0^\be\frac{2+\log(1+s\d^{-1})}{(1+\log(1+s\d^{-1}))^2(1+s\d^{-1})^2}\,f(s)\,ds\\[.2cm]
=&\int_0^\be\frac{2+\log(1+s)}{(1+\log(1+s))^2(1+s)^2}\,f(s\d)\,ds\\[.2cm]
\le&\,2\int_0^\be\frac{1}{(1+\log(1+s))(1+s)^2}\,f(s\d)\,ds.
\end{align*}
It remains to observe that
\begin{align*}
\int_0^e\frac{1}{(1+\log(1+s))(1+s)^2}\,f(s\d)\,ds
&\le f(e\d)\int_0^e\frac{1}{(1+\log(1+s))(1+s)^2}\,ds\\[.2cm]
&\le f(e\d)\int_0^\be\!\frac{ds}{(1+s)^2}=f(e\d)\le\const\int_e^\be\!\frac{f(s\d)\,ds}{s^2\log s}
\end{align*}
and
\bey
\int_e^\be\frac{1}{(1+\log(1+s))(1+s)^2}\,f(s\d)\,ds\le\int_e^\be\frac{f(s\d)\,ds}{s^2\log s}. \quad\bl
\eey

\begin{cor}
\label{lipconv}
Suppose that under the hypotheses of Theorem {\rm\ref{ostar}}, the function $f$ is bounded and has finite right derivative at $0$. Then
$$
\O_f(\d)\le \const a \d\,\log\Big(\log \frac M{a\d}\Big)\,\quad\mbox{for}\quad
\d\in\Big(0,\frac M{3a}\Big),
$$
where $a=f^\prime_+(0)$ and $M=\sup f$.
\end{cor}

\Pf Since $f(t)\le\min\{at,M\}$, $t>0$, the result follows from Theorem \ref{ostar} and the following obvious facts:
\begin{align*}
\int_e^{\frac M{a\d}}\frac{a\d ds}{s\log s}=a\d\log\left(\log\frac{M}{a\d}\right)
\quad\mbox{and}\quad
\int_{\frac M{a\d}}^\be\frac{M ds}{s^2\log s}\le\int_{\frac M{a\d}}^\be\frac{M ds}{s^2}=a\d.\quad\bl
\end{align*}

In \cite{AP2} we proved that if $f$ belongs to the H\"older class $\L_\a(\R)$, $0\le\a<1$,
then
\bay
\label{ato1}
\O_f(\d)\le\const(1-\a)^{-1}\|f\|_{\L_\a}\d^\a,\quad\d>0,
\ey
where
$$
\|f\|_{\L_\a}\df\sup_{x\ne y}\frac{|f(x)-f(y)|}{|x-y|}.
$$
The next result shows that if in addition to this $f$ satisfies the hypotheses of Theorem \ref{ostar}, then the factor $(1-\a)^{-1}$ on the right-hand side of \rf{ato1} can considerably be improved.

\begin{cor}
\label{holconv}
Suppose that under the hypotheses of Theorem {\em\ref{ostar}}, the function $f$ belongs to $\L_\a(\R)$, $0\le\a<1$.
Then
$$
\O_f(\d)\le \const\Big(\log\frac2{1-\a}\Big)\|f\|_{\L_\a}\d^{\a}
$$
for every $\d>0$.
\end{cor}

\Pf Indeed,
$$
\int_e^\be\frac{ds}{s^{2-\a}\log s}=\int_1^\be e^{(\a-1)t}\frac{dt}t
=\int_{1-\a}^\be \frac{e^{-t}\,dt}t\le\const\log\frac2{1-\a}.\quad\bl
$$

\medskip

{\bf Remark.} The function $x\mapsto1+\vk(x-1)$ satisfies the hypotheses of Corollary \ref{lipconv}
with $a=1$ and $M=2$, and Corollary \ref{lipconv} yields a sharp result in this case.
That means that Theorem \ref{ostar}
is also sharp in a sense.

\medskip

The following theorem is a symmetrized version of Theorem \ref{ostar}.

\begin{thm} Let $f$ be a continuous function on $\R$
such that $f$ is convex or concave on each of two rays $(-\be,0]$ and  $[0,\be)$.
Suppose that there exists a finite limit $\lim\limits_{|t|\to\be}t^{-1}f(t)\df a$.
Then
$$
\O_f(\d)\le a\d+ c\int_e^\be\frac{|f(\d s)-f(0)-\d as|+|f(-\d s)-f(0)+\d as|}{s^2\log s}\,ds,
$$
where $c$ as a numerical constant.
\end{thm}

\Pf It suffices to consider the case where $f(0)=a=0$.
We assume first that $f(t)=0$ for $t\le0$. To be definite, suppose that $f$ is
concave on $[0,\be)$.  Then $f$ is a nondecreasing function because $\lim\limits_{|t|\to\be}t^{-1}f(t)=0$, and so the result reduces to Theorem \ref{ostar}.
The case $f(t)=0$ for $t\ge0$ follows from the considered case with the help of
the change of variables $t\mapsto-t$.
It remains to observe that each function $f$ with $a=f(0)=0$
can be represented in the form $f=g+h$ in such way that
$g(t)=0$ for $t\le0$, $h(t)=0$ for $t\ge0$, and the cases of the function $g$ and $h$
have been treated above. $\bl$

\begin{thm}
\label{911}
Let $f$ be a nonnegative continuous function on $\R$ such that
$f(x)=0$ for all $x\le0$ and the function $x\mapsto x^{-1}f(x)$ is nonincreasing on $(0,\be)$.
Suppose that $\O_f(\d)<\be$ for $\d>0$. Then
$$
f(x)\le\const\frac x{\log\log x}
$$
for every $x\ge4$.
\end{thm}

\Pf By Theorem \ref{kme-},
$$
\O^\flat_f(1)\ge\frac12\|\dg_0 f\|_{\frak M_{[1,\be),(-\be,0]}}.
$$
Making the change of variables $y\mapsto-y$ we get
$$
\left\|\frac{f(x)}{x+y}\right\|_{\frak M_{[1,\be),[0,\be)}}\le2\O^\flat_f(1).
$$
Thus for every $a>1$
\bey
\left\|\frac{x}{x+y}\right\|_{\frak M_{[1,a],[1,a]}}
\le\max_{[1,a]}\left|\frac x{f(x)}\right|\cdot
\left\|\frac{f(x)}{x+y}\right\|_{\frak M_{[1,a],[1,a]}}\\
\le
\frac{a}{f(a)}\left\|\frac{f(x)}{x+y}\right\|_{\frak M_{[1,\be),[0,\be)}}
\le\frac{2a\O^\flat_f(1)}{f(a)}.
\eey
It remains to apply Theorem \ref{albe-}. $\bl$

\medskip

{\bf Remark.} Let $x_0>e$ and let $g_\a$ be a continuous function such that
$$
g_\a(x)=\left\{\begin{array}{ll}\dfrac x{\log^\a(\log x)},&\text {if}\,\,\,\,x\ge x_0>0,\\[.6cm]
0,&\text {if}\,\,\,\,x\le0.
\end{array}\right.
$$
Then $\O_{g_\a}(\d)<\be$ for $\a>1$. Indeed, in this case $g_\a$ coincides with a function satisfying Theorem \ref{ostar} outside a compact subset of $\R$. On the other hand,
$\O_{g_\a}(\d)=\be$ for $\a<1$. This follows from Theorem \ref{911}. Indeed, outside a compact subset of $\R$ the function $g_\a$ coincides with a function $f$, for which the function $x\mapsto x^{-1}f(x)$ is nonincreasing on $(0,\be)$.
The case $\a=1$ is an open problem.

\

\section{\bf Lower estimates for operator moduli of continuity}
\setcounter{equation}{0}
\label{lest}

\

Recall that it follows from \rf{modnep} that if $f$ is a function on $\R$ such that $\|f\|_{L^\be}\le1$, $\|f\|_{\Li}\le1$, then
$$
\O_f(\d)\le\const\,\d\left(1+\log\frac1\d\right),\quad\d\in(0,1].
$$
It is still unknown whether this estimate is sharp. In particular, the question whether one can replace the factor $\left(1+\log\frac1\d\right)$ on the right-hand side with $\left((1+\log\frac1\d\right)^s$
for some $s<1$ is still open.

In \S\,\ref{|x|} we established a lower estimate for the operator modulus of continuity of the function $x\mapsto|x|$ on finite intervals.

The main purpose of this section is to construct a $C^\be$ function $f$ on $\R$ such that
$\|f\|_{L^\be}\le1$, $\|f\|_{\Li}\le1$, and
$$
\O_f(\d)\ge\const\,\d\sqrt{\log\frac2\d},\quad\d\in(0,1].
$$

Let $\s>0$. Denote by $\mathscr E_\s$ the set of entire functions
of exponential type at most $\s$.

Let $F\in\mathscr E_\s\cap L^2(\R)$. Then
$$
F(z)=\sum_{n\in\Z}\frac{\sin(\s z-\pi n)}{\s z-\pi n}F\left(\frac{\pi n}{\s}\right),
$$
see, e.g., \cite{L}, Lect. 20.2, Th. 1.
Let $f\in\mathscr E_\s\cap L^\be(\R)$. Then
$$
f(z)\dfrac{\sin\big(\s(z-a)\big)}{\s(z-a)}\in\mathscr E_{2\s}\cap L^2(\R).
$$
Hence,
\begin{align*}
f(z)\frac{\sin(\s(z-a))}{\s(z-a)}&=\sum_{n\in\Z}\frac{\sin(2\s z-\pi n)}{2\s z-\pi n}\cdot
\frac{\sin\big(\s(\frac{\pi n}{2\s}-a)\big)}{\s\big(\frac{\pi n}{2\s}-a\big)}f\left(\frac{\pi n}{2\s}\right)\\[.2cm]
&=2\sum_{n\in\Z}\frac{\sin(2\s z-\pi n)\sin\big(\s a-\frac{\pi n}2\big)}{(2\s z-\pi n)(2\s a-\pi n)}
f\left(\frac{\pi n}{2\s}\right).
\end{align*}
Substituting $z=a$, we obtain
\begin{align}
\label{t}
f(z)&=2\sum_{n\in\Z}\frac{\sin(2\s z-\pi n)\sin(\s z-\frac{\pi n}2)}{(2\s z-\pi n)^2}
f\left(\frac{\pi n}{2\s}\right)\nonumber\\[.2cm]
&=\sum_{n\in\Z}\frac{\sin^2(\s z-\frac{\pi n}2)\cos(\s z-\frac{\pi n}2)}{(\s z-\frac{\pi n}2)^2}
f\left(\frac{\pi n}{2\s}\right)
\end{align}
for $f\in\mathscr E_\s\cap L^\be(\R)$.

Denote by $\mathscr E_\s\big(\C^2\big)$ the set of all entire functions $f$ on $\C^2$
such that the functions $z\mapsto f(z,\xi)$ and $z\mapsto f(\xi,z)$ belong to $\E_\s$ for every
$\xi\in\R$ (or, which is the same, for all $\xi\in\C$).
Equality \rf{t} implies the following identity:

\bay
\label{tt}
f(z,w)\!=\!\!\!\!
\sum_{(m,n)\in\Z^2}\!\!\!\!\frac{\sin^2(\s z\!-\!\frac{\pi m}2)\!\cos(\s z\!-\!\frac{\pi m}2)
\sin^2(\s w\!-\!\frac{\pi n}2)\cos(\s w\!-\!\frac{\pi n}2)}{(\s z-\frac{\pi m}2)^2(\s w-\frac{\pi n}2)^2}
f\!\left(\frac{\pi m}{2\s},\frac{\pi n}{2\s}\right)
\ey
for every $f\in\mathscr E_\s(\C^2)\cap L^\be(\R^2)$.

\begin{thm}
\label{mzr}
Let $\s>0$ and $\Phi\in\mathscr E_\s(\C^2)$.
Suppose that $\Phi(\frac{\pi m}{2\s}+\a,\frac{\pi n}{2\s}+\b)\in\frak M_{\Z,\Z}$
for some $\a,\b\in\R$.
Then $\Phi\in\frak M_{\R,\R}$ and
$$
\|\Phi(x,y)\|_{\frak M_{\R,\R}}
\le2\Big\|\Phi\Big(\frac{\pi m}{2\s}+\a,\frac{\pi n}{2\s}+\b\Big)\Big\|_{\frak M_{\Z,\Z}}.
$$
\end{thm}

\Pf Clearly, it suffices to consider the case when $\a=\b=0$, $\s=\pi/2$ and $\|\Phi(m,n)\|_{\frak M_{\Z,\Z}}=1$.
Then (see \cite{Pi}, Theorem 5.1) there exist two sequences $\{\var_m\}_{m\in\Z}$ and $\{\psi_n\}_{n\in\Z}$
of vectors in the closed unit ball of a Hilbert space $\mathcal H$ such that
$(\var_m,\psi_n)=\Phi(m,n)$.
Put
$$
g_x\df\frac4{\pi^2}\sum_{m\in\Z}\frac{\sin^2\big(\frac\pi2(x-m)\big)\cos\big(\frac\pi2(x-m)\big)}{(x-m)^2}
\f_m
$$
and
$$
h_y\df\frac4{\pi^2}\sum_{n\in\Z}\frac{\sin^2\big(\frac\pi2(y-n)\big)\cos\big(\frac\pi2(y-n)\big)}{(y-n)^2}
\psi_n.
$$
We have
\begin{align*}
\|g_{x}\|_{\mathcal H}&\le\frac4{\pi^2}\sum_{m\in\Z}\frac{\sin^2\big(\frac\pi2(x-m)\big)
|\cos\big(\frac\pi2(x-m)\big)|}{(x-m)^2}\\[.2cm]
&=\frac4{\pi^2}\sum_{n\in\Z}\frac{\sin^2\frac{\pi x}2\big|\cos\frac{\pi x}2\big|}{(x-2 n)^2}+
\frac4{\pi^2}\sum_{n\in\Z}\frac{\sin^2\big(\frac{\pi x}2-\frac\pi2\big)
\big|\cos\big(\frac{\pi x}2-\frac\pi2\big)\big|}{(x-2 n-1)^2}\\[.2cm]
&=\left|\cos\frac{\pi x}2\right|+\left|\sin\frac{\pi x}2\right|\le\sqrt2.
\end{align*}
In the same way, $\|h_y\|_{\mathcal H}\le \sqrt2$ for all $x\in\R$.
Clearly $|\Phi|\le1$ on $\Z^2$. The Cartwright theorem (see \cite{L}, Lecture 21, Theorem 4) implies
that $\Phi$ is bounded on $\R\times\Z$. Applying once more the Cartwright theorem, we find
that $\Phi\in L^\be(\R^2)$. Hence, we can apply formula \rf{tt} to the function $\Phi$, whence
$\Phi(x,y)=(g_x,h_y)$ for all $x, y\in\R$. It remains to observe that
by Theorem 5.1 in \cite{Pi},
$$
\|\Phi(x,y)\|_{\frak M_{\R,\R}}\le\sup_{x\in\R}\|g_x\|_{\mathcal H}\cdot
\sup_{y\in\R}\|h_y\|_{\mathcal H}\le2.\quad\bl
$$

\begin{thm} Let $f\in\mathscr E_{\s}$. Then
$$
\O_{f}^{\flat}(\d)\ge\frac\d{2}\left\|\frac{f(x)-f(y)}{x-y}\right\|_{\frak M_{\R,\R}}
$$
for every $\d\in\big(0,\frac1{2\s}\big]$.
\end{thm}

\Pf The general case easily reduces to the case $\s=\pi/4$.
By Theorem \ref{mzr}, we have
$$
\left\|\frac{f(x)-f(y)}{x-y}\right\|_{\frak M_{\R,\R}}\le2
\left\|\frac{f(2m+1)-f(2n)}{2m-2n+1}\right\|_{\frak M_{\Z,\Z}}\le2\|f\|_{{\rm OL}(\Z)}.
$$

Hence, by Theorem \ref{pi2},
$$
\O_{f}^{\flat}(\d)\ge\O_{f,\Z}^{\flat}(\d)=
\d \|f\|_{{\rm OL}(\Z)}\ge\frac\d2\left\|\frac{f(x)-f(y)}{x-y}\right\|_{\frak M_{\R,\R})}
$$
for $\d\in\big(0,\frac2\pi\big]$. $\bl$

\begin{thm}
\label{57}
Let $f\in\mathscr E_{\s}$. Then
$$
\O_{f}(\d)\ge\frac\d{4}\left\|\frac{f(x)-f(y)}{x-y}\right\|_{\frak M_{\R,\R}}
$$
for every $\d\in\big(0,\frac1{2\s}\big]$.
\end{thm}

\Pf It suffices to observe that $\O_{f}^{\flat}(\d)\le2\O_{f}(\d)$
by Theorem 10.2 in \cite{AP2}. $\bl$

\begin{thm}
\label{Log}
For every $\d\in(0,1]$, there exists an entire function $f\in\mathscr E_{1/\d}$
such that $\|f\|_{L^\be(\R)}\le1$, $\|f^\prime\|_{L^\be(\R)}\le1$ and
$\O_f(\d)\ge C\,\d\sqrt{\log\frac2\d}$, where $C$ is a positive numerical constant.
\end{thm}

We need some lemmata.

\begin{lem}
\label{trig}
For every positive integer $n$,
there exists a trigonometric polynomial $f$ of degree $n$ such that $\|f\|_{L^\be}\le1$,
$\|f^\prime\|_{L^\be}\le1$,
and
$$
\left\|\frac{f(x)-f(y)}{e^{{\rm i}x}-e^{{\rm i}y}}\right\|_{{\fM}_{[0,2\pi],[0,2\pi]}}
\ge c\,\sqrt{\log n}.
$$
\end{lem}

\Pf It follows from the results of
\cite{Pe1} that for every function $h$ in $C^1(\T)$,
\bay
\label{bes}
\left\|\frac{h(e^{{\rm i}x})-h(e^{{\rm i}y})}{e^{{\rm i}x}-e^{{\rm i}y}}\right\|_{{\fM}_{[0,2\pi],[0,2\pi]}}
\ge\const\|h\|_{B_1^1},
\ey
where $B_1^1$ is a Besov space (see \cite{Pe4} for the definition) of functions on $\T$.
Note that this result was deduced in \cite{Pe1} from the nuclearity criterion for Hankel
operators (see \cite{Pe0} and \cite{Pe4}, Ch. 6). It is easy to see from the definition
of $B_1^1(\T)$ (see e.g., \cite{Pe4}) that
\bay
\label{lacuna}
\|h\|_{B_1^1}\ge\const\sum_{j\ge0}2^j|\hat h(2^j)|.
\ey

It is well known (see, for example, \cite{Fo}) that for every positive integer $n$, there  exists an analytic polynomial $h$ such that
$$
h(0)=0,\quad\deg h=n,\quad\|h^\prime\|_{L^\be(\T)}=1,\quad\mbox{and}\quad
\sum\limits_{j\ge0}2^j|\hat h(2^j)|\ge d\,\sqrt{\log n},
$$
where $d$ is a positive numerical constant.
Then inequality \rf{bes} implies that
$$
\left\|\frac{h(e^{{\rm i}x})-h(e^{{\rm i}y})}{e^{{\rm i}x}-e^{{\rm i}y}}\right\|_{{\fM}_{[0,2\pi],[0,2\pi]}}
\ge\const\sqrt{\log n}.
$$
Put $f(x)\df h(e^{{\rm i}x})$. It remains to observe that
$\|f^\prime\|_{L^\be}=\|h^\prime\|_{L^\be(\T)}=1$ and
$\|f\|_{L^\be}=\|h\|_{L^\be(\T)}\le1$. $\bl$

\begin{lem}
\label{kva}
Let $n\in\Z$. Then
$$
\left\|\frac{x-y-2\pi n}{e^{{\rm i}x}-e^{{\rm i}y}}\right\|_{{\fM}_{J_1,J_2}}
\le\frac{3\sqrt2\,\pi}4
$$
for every intervals $J_1$ and $J_2$ with $J_1-J_2\subset\big[(2n-\frac32)\pi,(2n+\frac32)\pi\big]$.
\end{lem}

\Pf We can restrict ourselves to the case $n=0$.
We have
\begin{align*}
\left\|\frac{x-y}{e^{{\rm i}x}-e^{{\rm i}y}}\right\|_{{\fM}_{J_1,J_2}}
&=\left\|\frac{x-y}{e^{{\rm i}(x-y)}-1}\right\|_{{\fM}_{J_1,J_2}}\\[.2cm]
&\le\left\|\frac{t}{e^{{\rm i}t}-1}\right\|_{\widehat L^1([-\frac{3\pi}2,\frac{3\pi}2])}
=\left\|\frac{t}{2\sin(t/2)}\right\|_{\widehat L^1([-\frac{3\pi}2,\frac{3\pi}2])}.
\end{align*}
Consider the $3\pi$-periodic function $\xi$ such that $\xi(t)=\frac t{2\sin(t/2)}$
for $t\in[-\frac{3\pi}2,\frac{3\pi}2]$. We can expand $\xi$ in Fourier series
$$
\xi(t)=\sum_{n\in\Z}a_ne^{\frac23n{\rm i}t}.
$$
Note that $a_n=a_{-n}\in\R$ for all $n\in\Z$ because $\xi$ is even and real.
Moreover, $\xi$ is convex on $[-\frac{3\pi}2,\frac{3\pi}2]$. Hence, by Theorem 35 in \cite{HR}, $(-1)^na_n\ge0$ for all $n\in\Z$.
It follows that
$$
\left\|\frac{t}{2\sin(t/2)}\right\|_{\widehat L^1([-\frac{3\pi}2,\frac{3\pi}2])}
\le\sum_{n\in\Z}|a_n|=\xi\left(\frac{3\pi}2\right)=\frac{3\sqrt2\,\pi}4.\quad\bl
$$

\begin{cor}
\label{kvas}
Let $J_1=[\pi j,\pi j+\pi)]$ and $J_2=[\pi k -\frac\pi2,\pi k +\frac\pi2]$, where $j,k\in\Z$.
Then
$$
\left\|\frac{x-y-2\pi n}{e^{{\rm i}x}-e^{{\rm i}y}}\right\|_{{\fM}_{J_1,J_2}}\le\frac{3\sqrt2\,\pi}4
$$
for some $n\in\Z$.
\end{cor}

\Pf We have $J_1-J_2=[\pi(j-k)-\frac\pi2,\pi(j-k)+\frac{3\pi}2]$.
If $j-k$ is even, then we can apply Lemma \ref{kva} with $n=\frac12(j-k)$.
If $j-k$ is odd, then we can apply Lemma \ref{kva} with $n=\frac12(j-k+1)$. $\bl$

\begin{lem}
\label{kvad}
Let $g$ be a $2\pi$-periodic function in $C^1(\R)$. Then
$$
\left\|\frac{g(x)-g(y)}{e^{{\rm i}x}-e^{{\rm i}y}}\right\|_{{\fM}_{[0,2\pi],[0,2\pi]}}
\le3\sqrt2\,\pi\left\|\frac{g(x)-g(y)}{x - y}\right\|_{{\fM}_{\R,\R}}.
$$
\end{lem}
\Pf Note that
$$
\left\|\frac{g(x)-g(y)}{x - y}\right\|_{{\fM}_{\R,\R}}=\left\|\frac{g(x)-g(y)}{x - y-2\pi n}\right\|_{{\fM}_{\R,\R}}
$$
for all $n\in\Z$ and
$$
\left\|\frac{g(x)-g(y)}{e^{{\rm i}x}-e^{{\rm i}y}}\right\|_{{\fM}_{[0,2\pi],[0,2\pi]}}=
\left\|\frac{g(x)-g(y)}{e^{{\rm i}x}-e^{{\rm i}y}}\right\|_{{\fM}_{[0,2\pi],[-\frac\pi2,\frac{3\pi}2]}}.
$$
Now we can represent the square $[0,2\pi]\times[-\frac\pi2,\frac{3\pi}2]$ as the union of four squares with sides of length $\pi$, each of which satisfies the hypotheses of Corollary \ref{kvas}. $\bl$

\medskip

{\bf Proof of Theorem \ref{Log}.} It suffices to consider the case when $\d\in\big(0,\frac12\big]$.
Then $\d\in\big[\frac1{n+1},\frac1n\big]$ for an integer $n\ge2$. By Lemma \ref{trig},
there exists a trigonometric polynomial $f$ of degree $n$ such that $\|f\|_{L^\be}\le1$,
$\|f^\prime\|_{L^\be}\le1$
and
$$
\left\|\frac{f(x)-f(y)}{e^{{\rm i}x}-e^{{\rm i}y}}\right\|_{{\fM}_{[0,2\pi],[0,2\pi]}}
\ge c\sqrt{\log n}.
$$
Hence,
$$
\left\|\frac{f(x)-f(y)}{x-y}\right\|_{{\fM}_{\R,\R}}
\ge c\sqrt{\log n}
$$
by Lemma \ref{kvad}.
Clearly,
$g\in\E_n\subset\E_{1/\d}$.
Applying Theorem \ref{57},
we obtain
$$
\O_f(t)\ge\const\sqrt{\log n}\,\, t,\quad 0<t\le\frac1{2n}.
$$
Hence,
$$
\O_f(\d)\ge\O_f\left(\frac1{2n}\right)\ge C_0\frac{\sqrt{\log n}}{n}\ge C\d\sqrt{\log\left(\frac2\d\right)}
$$
for some positive numbers $C_0$ and $C$. $\bl$

\begin{thm}
\label{nizko}
There exist a positive number $c$ and a function $f\in C^\be(\R)$ such that $\|f\|_{L^\be}\le1$,
$\|f^\prime\|_{L^\be}\le1$, and $\O_f(\d)\ge c\,\d\sqrt{\log\frac2\d}$ for every $\d\in(0,1]$.
\end{thm}

\Pf Applying Theorem \ref{Log} for $\d=2^{-n}$, we can construct a sequence functions
$\{f_n\}_{n\ge1}$ and two sequences of bounded self-adjoint operators $\{A_n\}_{n\ge1}$
and $\{B_n\}_{n\ge1}$ such that $\|f_n\|_{L^\be}\le1$, $\|f_n^\prime\|_{L^\be}\le1$, $\|A_n-B_n\|\le 2^{-n}$ and
$\|f_n(A_n)-f_n(B_n)\|\ge C \sqrt n\,\, 2^{-n}$ for all $n\ge1$. Denote by $\D_n$ the convex hull
of $\s(A_n)\cup\s(B_n)$. Using the translations
$f_n\mapsto f_n(x-a_n)$, $A_n\mapsto A_n+a_nI$, $B_n\mapsto B_n+a_nI$
and $\D_n\mapsto a_n+\D_n$ for a suitable sequence $\{a_n\}_{n=1}^\be$ in $\R$, we can achieve the condition that the intervals $\D_n$ are
disjoint and $\dist(\D_m,\D_n)>2$ for $m\not=n$.
We can construct a function $f\in C^\be(\R)$ such that $\|f\|_{L^\be}\le1$, $\|f^\prime\|_{L^\be}\le1$ and
$f\big|\D_n=f_n\big|\D_n$ for all $n\ge1$.
Clearly, $\O_f(2^{-n})\ge C \sqrt n\,\, 2^{-n}$ for all $n\ge1$ and some positive $C$
which easily implies the result. $\bl$

To obtain the lower estimate in Theorem \ref{nizko}, we used the inequality
\bay
\label{slabaya}
\left\|\frac{f(e^{{\rm i}x})-f(e^{{\rm i}y})}{e^{{\rm i}x}-e^{{\rm i}y}}\right\|_{{\fM}_{[0,2\pi],[0,2\pi]}}
\ge\const\sum_{j\ge0}2^j|\hat f(2^j)|,
\ey
which in turn implies that there exists a positive number $C$ such that for every positive integer $n$ there exists a polynomial $f$ of degree $n$ such that
\bay
\label{kor}
\left\|\frac{f(e^{{\rm i}x})-f(e^{{\rm i}y})}{e^{{\rm i}x}-e^{{\rm i}y}}\right\|_{{\fM}_{[0,2\pi],[0,2\pi]}}
\ge C\sqrt{\log n}\,\|f\|_{\Li}.
\ey

We do not know whether Theorem \ref{nizko} can be improved. It would certainly be natural to try to improve \rf{kor}. The best known lower estimate for the norm of divided differences in the space of Schur multipliers was obtained in \cite{Pe1}. To state it, we need some definitions.

Let $f\in L^1(\T)$. Denote by $\cp f$ the Poisson integral of $f$,
$$
(\cp) f(z)\df\int_\T\frac{(1-|z|^2)f(\z)}{|z-\z|^2}\,d\m(\z),\quad z\in\dd,
$$
where $\m$ is normalized Lebesgue measure on $\T$.

For $t\in\R$ and $\d\in(0,1)$, we define the Carleson domain $Q(t,\d)$ by
$$
Q(t,\d)\df\{re^{{\rm i}s}:0<1-r<h,|s-t|<\d\}.
$$
A positive  Borel measure on $\mu$ on $\dd$ is said to be a {\it Carleson measure}
if
$$
\cc(\mu)\df\mu(\dd)+\sup\big\{\d^{-1}\mu(Q(t,\d)):~t\in\R,~\d\in(0,1)\big\}<\be.
$$
If $\psi$ is a nonnegative measurable function on $\dd$, we put
$$
\cc(\psi)\df\cc(\mu),\quad\mbox{where}\quad d\mu\df\psi\,d\m_2.
$$
Here $\m_2$ is planar Lebesgue measure.

It follows from results of \cite{Pe1} (see also \cite{Pe}) that
\bay
\label{sil'naya}
\left\|\frac{f(e^{{\rm i}x})-f(e^{{\rm i}y})}{e^{{\rm i}x}-e^{{\rm i}y}}\right\|_{{\fM}_{[0,2\pi],[0,2\pi]}}\ge\const\|f\|_{\cL},
\ey
where
$$
\|f\|_\cL\df\cc\big(\big\|{\rm Hess}(\cp f)\big\|\big),
$$
where for a function $\f$ of class $C^2$, its Hessian ${\rm Hess}(\f)$ is the matrix of its second order partial derivatives.

It turns out, however, that for a trigonometric polynomial $f$ of degree $n$,
\bay
\label{kuzya}
\|f\|_\cL\le\const\sqrt{\log(1+n)}\|f\|_\Li,
\ey
and so even if instead of inequality \rf{slabaya} we use inequality \rf{sil'naya}, we cannot improve Theorem \ref{nizko}.

Inequality \rf{kuzya} is an immediate consequence of the following fact:

\begin{thm}
\label{Carl}
For a trigonometric polynomial $f$ of degree $n$, $n\ge2$, the following inequality holds:
$$
\cc\big(\big|\nabla(\cp f)\big|\big)\le\const\sqrt{\log n}\,\|f\|_{L^\be}.
$$
\end{thm}

We are going to use the well-known fact that a function $f$ in $L^1(\T)$ belongs to the space ${\rm BMO}(\T)$ if
and only if the measure $\mu$ defined by $d\,\mu=|\nabla(\cp f)|^2(1-|z|)\,d\m_2$ is a Carleson measure.
We refer to \cite{Ga} for Carleson measures and the space ${\rm BMO}$.

\medskip

{\bf Proof of Theorem \ref{Carl}.} Suppose that $\|f\|_{L^\infty}=1$. We have to prove that
\bay
\label{logn}
\int_{Q(t,\d)}
|\nabla(\cp f)|\,dxdy\le\const\,\d\sqrt{\log n}\,.
\ey
Note that $|\nabla(\mathscr P f)|\le2n$ by Bernstein's inequality. Hence,
\begin{align*}
\int_{\{1-n^{-1}<|\z|<1\}\cap Q(t,\d)}
|\nabla(\cp f)|\,d\m_2
&\le 2n\m_2\big(\{\z:~1-n^{-1}<|\z|<1\}\cap Q(t,\d)\big)\\[.2cm]
&=2n\d(1-(1-n^{-1})^2)\le4\d.
\end{align*}
This proves \rf{logn} in the case $\d\ge1- n^{-1}$.
In the case  $\d<1-n^{-1}$ it remains to estimate the integral over the set
$\{\z:~|\z|<1-n^{-1}\}\cap Q(t,\d)$. Note that
$\|f\|_{\rm BMO}\le\const\|f\|_{L^\be}$. Hence, there exists a constant $C$
such that
$$
\int_{Q(t,\d)}|\nabla(\mathscr P f)|^2(1-|\z|)\,d\m_2(\z)\le C\d.
$$
Thus
\begin{align*}
\int_{\{|\z|<1-n^{-1}\}\cap Q(t,\d)}&
|\nabla(\cp f)|\,d\m_2\\[.2cm]
&\hspace*{-2.55cm}\le\left(\int_{Q(t,\d)}
|\nabla(\cp f)|^2(1-|\z|)\,d\m_2(\z)\right)^{1/2}\!\!
\left(\int_{\{|\z|<1-n^{-1}\}\cap Q(t,\d)}
(1-|\z|)^{-1}\,d\m_2(\z)\right)^{1/2}\\[.2cm]
&\hspace*{-2.55cm}\le\const\d(\log(n\d))^{1/2}
\le \const\d(\log n)^{1/2}.\quad\bl
\end{align*}

\

\section{\bf Lower estimates in the case of unitary operators}
\setcounter{equation}{0}
\label{lest+}

\

The purpose of this section is to obtain lower estimates for the operator modulus of continuity for functions on the unit circle.

We define an operator modulus of continuity of a continuous function $f$ on $\T$ by
$$
\O_f(\d)\df\sup\big\{\|f(U)-f(V)\|:~U~\mbox{and}~V~\mbox{are unitary},~ \|U-V\|\le\d\big\}.
$$
As in the case of self-adjoint operators (see \cite{AP2}),
one can prove that
$$
\|f(U)R-Rf(V)\|\le2\O_f(\|UR-RV\|)
$$
for every unitary operators $U$, $V$ and an operator $R$
of norm $1$.
We define the space ${\rm OL}(\T)$ as the set of $f\in C(\T)$ such that
$$
\|f\|_{{\rm OL}(\T)}\df\sup_{\d>0}\d^{-1}\O_f(\d)<\be.
$$
Given a closed subset $\fF$ of $\T$, we can also introduce the operator modulus of continuity $\O_{f,\fF}$ and
define the space ${\rm OL}(\fF)$ of operator Lipschitz functions on $\fF$.

For closed subsets $\fF_1$ and $\fF_2$ of $\T$, the space $\fM_{\fF_1,\fF_2}$
of Schur multipliers can be defined by analogy with the self-adjoint case.
Note that the analogues of \rf{ravvo} and \rf{nva} for functions on closed subsets of $\T$
can be proved as in \S\,\ref{OL}.

Let $f\in C(\T)$. We put $f_\spadesuit(t)\df f(e^{{\rm i}t})$.
It is clear that $\O_{f_\spadesuit}\le\O_f$. Hence,
$\|f_\spadesuit\|_{{\rm OL}(\R)}\le\|f\|_{{\rm OL}(\T)}$.
Lemma \ref{kvad} implies that $\|f\|_{{\rm OL}(\T)}\le3\sqrt2\,\pi\|f_\spadesuit\|_{{\rm OL}(\R)}$.
One can prove that $\O_f\le\const\O_{f_\spadesuit}$.

Recall that it follows from results of \cite{Pe1} that for $f\in C(\T)$,
$$
\|f\|_{{\rm OL}(\T)}\ge\const\|f\|_{B^1_1};
$$
actually we used this estimate in \S\,\ref{lest}, see inequality \rf{bes}.

We would like to remind also that
for each positive integer $n$, there exists an analytic polynomial
$f$ such that $\deg f=n$, $\|f^\prime\|_{L^\be(\T)}=1$, and $\|f\|_{{\rm OL}(\T)}\ge\const\sqrt{\log n}$; see Lemma \ref{trig}.

Put
$$
\fd_n(z)\df\frac1n\frac{z^n-1}{z-1}=\frac1n\sum_{k=0}^{n-1}z^k.
$$
It is easy to see that
$$
\fd_n(\z z^{-1})=z^{1-n}\frac{z^n-\z^n}{n(z-\z)}=z^{1-n}\z^{n-1}\fd_n(z\z^{-1}).
$$

Denote by $\T_n$ the set of $n$th roots of 1, i.e., $\T_n\df\{\z\in\T:~\z^n=1\}$.

Let $f$ be an analytic polynomial of degree less that $n$.
Then
$$
f(z)=\sum_{\z\in\tau\T_n}f(\z)\fd_n(z\z^{-1})\quad\mbox{for every}\quad\t\in\T.
$$
If $f$ is a trigonometric polynomial and $\deg f\le n$,
then for every $\xi\in\T$, the function $z^nf(z)\fd_{2n}(z\xi^{-1})$ is an analytic polynomial of degree less than $4n$.
Hence,
$$
z^nf(z)\fd_{2n}(z\xi^{-1})=\sum_{\z\in\tau\T_{4n}}f(\z)\fd_{2n}(\z\xi^{-1})\fd_{4n}(z\z^{-1}).
$$
Substituting $\xi=z$ we get
\bay
\label{t+}
f(z)=z^{-n}\sum_{\z\in\t\T_{4n}}f(\z)\fd_{2n}(\z z^{-1})\fd_{4n}(z\z^{-1})
=\sum_{\z\in\t\T_{4n}}f(\z)F_n(z,\z)
\ey
for every $\t\in\T$, where
$$
F_n(z,\z)\df z^{1-3n}\z^{1-4n}\frac{(z^{2n}-\z^{2n})(z^{4n}-\z^{4n})}{8n^2(z-\z)^2}.
$$

Denote by $\mP_n\big(\T^2\big)$ the set of all trigonometric polynomial $f$ on $\T^2$
such that the functions $z\mapsto f(z,\xi)$ and $z\mapsto f(\xi,z)$ are trigonometric
polynomials on $\T$ of degree at most $n$  for every
$\xi\in\T$.
Equality \rf{t+} implies the following identity:
\bay
\label{tt+}
f(z,w)
=\sum_{\z\in\t_1\T_{4n}}\,\,\sum_{\xi\in\t_2\T_{4n}}f(\z,\xi)
F_n(z,\z)F_n(w,\xi)
\ey
for every $f\in\mP_n\big(\T^2\big)$ and for arbitrary $\t_1$ and $\t_2$ in $\T$.

\begin{thm}
\label{mzr+}
Let $\Phi\in\mP_n\big(\T^2\big)$. Then
$$
\|\Phi\|_{\fM_{\T,\T}}\le2\|\Phi\|_{{\fM}_{\t_1\T_{4n},\t_2\T_{4n}}}
$$
for every $\t_1,\t_2\in\T$.
\end{thm}
\Pf Clearly, it suffices to consider the case when $\t_1=\t_2=1$.
Then (see \cite{Pi}, Theorem 5.1) there exist two sequences $\{\var_\z\}_{\z\in\T_{4n}}$ and $\{\psi_\xi\}_{\xi\in\T_{4n}}$
of vectors in the closed unit ball of a Hilbert space $\mathcal H$ such that
$(\var_{\z},\psi_{\xi})=\Phi(\z,\xi)$.
Put
$$
g_z\df\sum_{\z\in\T_{4n}}F_n(z,\z)\f_\z\quad\text{and}\quad h_w\df\sum_{\xi\in\T_{4n}}F_n(w,\xi)\psi_\xi.
$$
Taking into account that for $z\in\T$,
\begin{align*}
\frac1{2n}\sum_{\z\in\T_{2n}}\left|\frac{z^{2n}-\z^{2n}}{z-\z}\right|^2
&=\frac1{2n}\sum_{\z\in\T_{4n}\setminus\T_{2n}}\left|\frac{z^{2n}-\z^{2n}}{z-\z}\right|^2\\[.2cm]
&=\int_\T\left|\frac{z^{2n}-\z^{2n}}{z-\z}\right|^2\,d\m(\z)=2n,
\end{align*}
we obtain
\begin{align*}
\|g_{z}\|_{\mathcal H}
&\le\sum_{\z\in\T_{4n}}|F_n(z,\z)|\\[.2cm]
&\le\frac{|z^{2n}+1|}{8n^2}\sum_{\z\in\T_{2n}}\left|\frac{z^{2n}-\z^{2n}}{z-\z}\right|^2
+\frac{|z^{2n}-1|}{8n^2}\sum_{\z\in\T_{4n}\setminus\T_{2n}}\left|\frac{z^{2n}-\z^{2n}}{z-\z}\right|^2\\[.2cm]
&=
\frac{|z^{2n}+1|+|z^{2n}-1|}2\le\sqrt2.
\end{align*}
In the same way, $\|h_w\|_{\mathcal H}\le \sqrt2$ for every $w\in\T$.
By \rf{tt+}, we have
$\Phi(z,w)=(g_z,h_w)$ for all $z, w\in\T$. It remains to observe that
by Theorem 5.1 in \cite{Pi},
$$
\|\Phi(z,w)\|_{\frak M_{\T,\T}}\le\sup_{z\in\T}\|g_z\|_{\mathcal H}\cdot
\sup_{w\in\T}\|h_w\|_{\mathcal H}\le2.\quad\bl
$$

We need the following version of Theorem \ref{pi2}:

\begin{thm}
\label{pi2+}
Let $f$ be a function on $\T_n$. Then
$$
\O_{f,\T_n}^{\flat}(\d)=\d\|f\|_{{\rm OL}(\T_n)}
$$
for every $\d\in(0,\frac4n]$.
\end{thm}

To prove Theorem \ref{pi2+}, we need a lemma. Put
$$
\l(z)\df\left\{\begin{array}{ll}z^{-1},&\text {if}\,\,\,z\in\C,\,\,z\ne0,\\[.2cm]
0,&\text {if}\,\,\,z=0.
\end{array}\right.
$$

\begin{lem}
\label{Hmn+}
Let $n$ be a positive integer.
Then
$$
\|\l(z-w)\|_{\fM_{\T_n,\T_n}}=\left\{\begin{array}{ll}\frac n4,&\text {if}\,\,\,\,n\,\,\,\text{is even},\\[.2cm]
\frac{n^2-1}{4n},&\text {if}\,\,\,\,n\,\,\,\text{is odd}.
\end{array}\right.
$$
\end{lem}
\Pf It is easy to verify that
$$
\sum_{k=1}^n\left(k-\frac{n+1}2\right)z^k=\frac{nz^n}{z-1}-\frac{z^n-1}{(z-1)^2}
-\frac{n+1}2 z\frac{z^n-1}{z-1}=n\l(z-1)
$$
for $z\in\T_n$. Hence,
\bay
\label{ft}
\l(z-w)=w^{-1}\l(zw^{-1}-1)=\frac1n\sum_{k=1}^n\left(k-\frac{n+1}2\right)z^kw^{-k-1}.
\ey
Thus
$$
\|\l(z-w)\|_{\fM_{\T_n,\T_n}}\le\frac1n\sum_{k=1}^n\left|k-\frac{n+1}2\right|=
\left\{\begin{array}{ll}\frac n4,&\text {if}\,\,\,\,n\,\,\,\text{is even},\\[.2cm]
\frac{n^2-1}{4n},&\text {if}\,\,\,\,n\,\,\,\text{is odd}.
\end{array}\right.
$$
The opposite inequality is also true. It can be deduced from the observation that equality
\rf{ft} means that that the function $\l(z-1)$ on the group $\T_n$ is the Fourier transform of the $n$-periodic
sequence $\{a_k\}_{k\in\Z}$ defined by $a_k=k-\frac{n+1}2$ for $k=1,2,\dots,n$.
Here we  identify the group dual to $\T_n$ with the group $\Z/n\Z$.
We omit details because we need only the upper estimate. $\bl$

\medskip

{\bf Proof of Theorem \ref{pi2+}.}  The inequality
$$
\O_{f,\T_n}^\flat(\d)\le\d\|f\|_{{\rm OL}(\T_n)},\quad\d>0,
$$
is a consequence of a unitary version of Theorem \ref{Olip}, which can be proved in the same way as the self-adjoint version, see also Theorem 4.13 in \cite{AP4}.

Let us prove the opposite inequality for $\d\in\big(0,\frac4n\big]$.
Fix $\e>0$. There exists a unitary operator $U$ and bounded operator $R$ such that
$\|UR-RU\|=1$, $\s(U)\subset\T_n$, and $\|f(U)R-Rf(U)\|\ge\|f\|_{{\rm OL}(\T_n)}-\e$.
Put
$$
R_U\df\sum_{\z,\xi\in\T_n,\,\z\not=\xi}E_U(\{\z\})RE_U(\{\xi\})=R-\sum_{\z\in\T_n}E_U(\{\z\})RE_U(\{\z\}).
$$
Clearly, $UR-RU=UR_U-R_UU$ and $f(U)R-Rf(U)=f(U)R_U-R_Uf(U)$.
Thus we may assume that $R=R_U$.
Note that
$$
UR-RU=\sum_{\z,\xi\in\T_n,\,\z\not=\xi}(\z-\xi)E_U(\{\z\})RE_U(\{\xi\}).
$$
Since
$$
R=R_U=\sum_{\z,\xi\in\T_n,\,\z\not=\xi}(\z-\xi)\l(\z-\xi)E_U(\{\z\})RE_U(\{\xi\}),
$$
we have $R=H_n\star(UR-RU)$, where $H_n(\z,\xi)=\l(\z-\xi)$, where $\z,\xi\in\T_n$.
Thus by Lemma \ref{Hmn+},
$$
\|R\|\le\|H_n\|_{\frak M_{\T_n,\T_n}}\|UR-RU\|=\|H_n\|_{\frak M_{\T_n,\T_n}}\le\frac n4.
$$

Let $\d\in\big(0,\frac4n\big]$. Then $\|U(\d R)-(\d R)U\|=\d$ and $\|\d R\|\le1$.
Hence,
$$
\O_{f,\T_n}^\flat(\d)\ge\d\|f(U)R-Rf(U)\|\ge\d\big(\|f\|_{{\rm OL}(\T_n)}-\e\big).
$$
Passing to the limit as $\e\to0$,
we obtain the desired result. $\bl$

\begin{thm}
Let $f$ be a trigonometric polynomial of degree $n\ge1$.
Then
$$
\O_{f,\T}^{\flat}(\d)\ge\frac\d2\|f\|_{{\rm OL}(\T)}
$$
for $\d\in(0,\frac1n]$.
\end{thm}

\Pf Applying Theorems \ref{mzr+} and \ref{pi2+}, we obtain
\begin{align*}
\left\|\frac{f(z)-f(w)}{z-w}\right\|_{\fM_{\T,\T}}
&\le2
\left\|\frac{f(z)-f(w)}{z-w}\right\|_{\fM_{\T_{4n},\T_{4n}}}\\[.2cm]
&=2\d^{-1}\O_{f,\T_{4n}}^{\flat}(\d)\le2\d^{-1}\O_{f,\T}^{\flat}(\d)
\end{align*}
for $\d\in(0,\frac1n]$. $\bl$

\begin{thm}
\label{polN}
Let $f$ be a trigonometric polynomial of degree $n\ge1$.
Then
$$
\O_{f,\T}(\d)\ge\frac\d4\|f\|_{{\rm OL}(\T)}
$$
for $\d\in(0,\frac1n]$.
\end{thm}
\Pf It suffices to observe that $\O_{f,\T}^{\flat}(\d)\le2\O_{f,\T}(\d)$. $\bl$

\begin{thm}
\label{besN}
Let $f\in C(\T)$. Then
$$
\O_f(2^{-n})\ge C\,2^{-n}\sum_{k=0}^{n-1}2^k\left(\big|\widehat f(2^k)\big|+\big|\widehat f(-2^k)\big|\right),
$$
where $C$ is a positive constant.
\end{thm}
\Pf Applying the convolution with the de la Vall\'ee Poussin kernel, we can find
an analytic polynomial $f_n$ such that $\deg f_n<2^{n}$,
$\widehat f_n(k)=\widehat f (k)$ for $k\le 2^{n-1}$ and $\O_{f_n}\le3\O_f$.
Applying inequalities \rf{bes} and \rf{lacuna}, we obtain
$$
\|f_n\|_{{\rm OL}(\T)}\ge\const\sum_{k=0}^{n-1}2^k\big(|\widehat f(2^k)|+|\widehat f(-2^k)|\big).
$$
It remains to apply Theorem \ref{polN} for $\d=2^{-n}$. $\bl$

In the following theorem we use the notation $C_A$ for the disk-algebra:
$$
C_A\df\big\{f\in C(\T):~\hat f(n)=0~\mbox{for}~n<0\big\}.
$$

\begin{thm}
\label{log2}
Let $\o:(0,2]\to\R$ be a positive continuous function. Suppose that
$\o(2t)\le \const\o(t)$, the function $t\mapsto t^{-1}(\log\frac4t)^{-1}\o(t)$ is nondecreasing, and
\bay
\label{oint}
\int_0^2\frac{\o^2(t)\,dt}{t^3\log^2\frac4t}<\be.
\ey
Then there exists a function $f\in C_A$ such that
$f^\prime\in C_A$ and $\O_f(\d)\ge\o(\d)$ for all $\d\in(0,2]$.
\end{thm}

\Pf Note that the inequality $\O_f(\d)\ge\o(\d)$ for $\d=2^{-n}$ implies that
$\O_f(\d)\ge\const\o(\d)$ for all $\d\in(0,2]$. Thus it
suffices to obtain the desired estimate for $\d=2^{-n}$. Taking Theorem \ref{besN} into account, we can reduce the result to the problem
to construct a function $g\in C_A$ such that
$$
a_n\df\frac{2^n\o(2^{-n})}n\le\frac1n\sum_{k=0}^{n-1}\big|\widehat g(2^k)\big|
$$
for all nonnegative integer $n$.

Indeed, in this case the function $f$ defined by
$$
f(z)=\int_0^z\frac{g(\z)-g(0)}\z\,d\z
$$
satisfies the inequality
$$
a_n\le\frac1n\sum_{k=0}^{n-1}2^k\big|\widehat f(2^k)\big|.
$$
Condition \rf{oint} implies that $\{a_n\}_{n\ge0}\in\ell^2$.
Moreover, $\{a_n\}_{n\ge0}$ is a nonincreasing sequence because the function
$t\mapsto t^{-1}(\log\frac4t)^{-1}\o(t)$ is nondecreasing.

We can find a function $g\in C_A$ such that $\widehat g(2^k)=a_k$ for all $k\ge0$,
see, for example, \cite{Fo}. Then
$$
\frac1n\sum_{k=0}^{n-1}\big|\widehat g(2^k)\big|=\frac1n\sum_{k=0}^{n-1}a_k\ge a_{n-1}\ge a_n.\quad\bl
$$

\medskip

{\bf Remark.} Theorem \ref{log2} remains valid if we replace the assumption
that the function $t\mapsto t^{-1}(\log\frac4t)^{-1}\o(t)$ is nondecreasing with the assumption that
there exists a positive constant $C$ such that
$$
\frac{\o(t)}{t\log\frac4t}
\le C \frac{\o(s)}{s\log\frac4s},\quad\mbox{whenever}\quad0<t<s\le2.
$$

\

\section{\bf Self-adjoint operators with finite spectrum.\\ Estimates in terms of the $\bs{\e}$-entropy of the spectrum}
\setcounter{equation}{0}
\label{safinsp}

\

In this section we obtain sharp estimates of the quasicommutator norms \lb$\|f(A)R-Rf(B)\|$ in the case when $A$ has finite spectrum. This allows us to obtain sharp estimates of the operator Lipschitz norm in terms of the Lipschitz norm in the case of operators on finite-dimensional spaces in terms of the dimension.

Moreover, we obtain a more general result (see Theorem \ref{qcom}) in terms of $\e$-entropy of the spectrum of $A$, where $\e=\|AR-RA\|$. This leads to an improvement of inequality \rf{modnep}.

Note that the results of this section improve some results of \cite{F2} and \cite{F3}.

Let $\fF$ be a closed subset of $\R$. Denote by ${\rm Lip}(\fF)$ the set
of Lipschitz functions on $\fF$. Put
$$
\|f\|_{{\rm Lip}(\fF)}\df\inf\big\{C>0:|f(x)-f(y)|\le C|x-y|\quad\forall x,\,y\in \fF\big\}.
$$

Let $\{s_j(T)\}_{j=0}^\be$ be the sequence of singular values of
a bounded operator $T$. We use the notation
$\bS_\o$ for the Matsaev ideal,
$$
\bS_\o\df\big\{T:\|T\|_{\bS_\o}\df\sum_{j=0}^\be(1+j)^{-1}s_j(T)<\be\big\}.
$$

We need the following statement which is contained implicitly in \cite{NP}.

\begin{thm}
\label{NP}
Let $f$ be a Lipschitz function on a closed subset $\fF$ of $\R$.
Then for every nonempty finite subset $\L$ in $\fF$,
$$
\|\dg_0 f\|_{{\frak M}_{\L,\fF}}\le C\big(1+\log({\rm card}(\L))\big)\|f\|_{{\rm Lip}(\fF)},
$$
 where $C$ is a numerical constant.
\end{thm}

\Pf Let $k\in L^2(\mu\otimes\nu)$, where $\mu$ and $\nu$ are Borel measures on $\L$ and $\fF$.
Clearly, $\rank\mI_k^{\mu,\nu}\le\card(\L)$.
Hence, $\|\mI_k^{\mu,\nu}\|_{\bS_\o}\le\big((1+\log(\card(\L))\big)\|\mI_k^{\mu,\nu}\|$.
Now Theorem 2.3 in \cite{NP} implies that
$$
\big\|\mI_{k\dg_0 f}^{\mu,\nu}\big\|\le\const\big((1+\log(\card(\L))\big)\big\|\mI_k^{\mu,\nu}\big\|\cdot\|f\|_{{\rm Lip}(\fF)}.
\quad\bl
$$

\begin{thm}
\label{skon}
Let $A$ and $B$ be self-adjoint operators. Suppose that $\s(A)$ is finite.
Then
$$
\|f(A)R-Rf(B)\|\le C\big(1+\log(\card(\s(A)))\big)\|f\|_{{\rm Lip}(\s(A)\cup\s(B))}\|AR-RB\|
$$
for all bounded operators $R$ and $f\in {\rm Lip}\big(\s(A)\cup\s(B)\big)$, where $C$ is a numerical constant.
\end{thm}

\Pf The result follows from Theorem \ref{NP} if we take into account the following
generalizations of \rf{doi} and \rf{doie} (see \cite{BS4}):
$$
f(A)R-Rf(B)=\iint\limits_{\s(A)\times\s(B)}(\dg_0f\big)(x,y)\,dE_A(x)(AR-RB)\,dE_B(y)
$$
and
$$
\left\|\,\,~\iint\limits_{\s(A)\times\s(B)}(\dg_0f\big)(x,y)\,dE_A(x)(AR-RB)\,dE_B(y)\right\|\le
\|\dg_0f\|_{{\fM}(\s(A)\times\s(B))}\|AR-RB\|
$$
which proves the result. $\bl$

\begin{cor}
\label{Cn}
Let $A$, $B$ be self-adjoint operators and let $R$ be a linear operator on $\C^n$.
Then
\bay
\label{ranee}
\|f(A)R-Rf(B)\|\le C(1+\log n)\|f\|_{{\rm Lip}(\s(A)\cup\s(B))}\|AR-RB\|
\ey
for every function $f$ on $\s(A)\cup\s(B)$, where $C$ is a numerical constant.
\end{cor}

\medskip

{\bf Remark 1.}
Note that in the special case $f(t)=|t|$ inequality \rf{ranee} is well-known, see, e.g., \cite{Da}.  This special case also follows from Matsaev's theorem, see \cite{GK}, Ch. III, Th. 4.2 (see also \cite{Go} where a finite dimensional improvement of Matsaev's theorem was obtained).

\medskip

{\bf Remark 2.}
We also would like to note that inequality \rf{ranee} is sharp. Indeed, it follows immediately from Lemma 15 of \cite{Da} that for each positive integer
$n$ there exist $n\times n$ self-adjoint matrices $A$ and $R$ such that
\bay
\label{MG}
\big\|\,|A|R-R|A|\,\big\|\ge\const\log(1+n)\|AR-RA\|\quad\text{and}\quad AR-RA\ne\0.
\ey
We also refer the reader to \cite{Mc} where inequality \rf{MG} is essentially contained.
Moreover, \rf{MG} can be deduced from the results of Matsaev and Gohberg mentioned above.

\medskip

The following result is a special case of Corollary \ref{Cn} that corresponds to $R=I$.

\begin{thm}
\label{Cnr}
Let $A$, $B$ be self-adjoint operators on $\C^n$.
Then
$$
\|f(A)-f(B)\|\le C(1+\log n)\|f\|_{{\rm Lip}(\s(A)\cup\s(B))}\|A-B\|
$$
for every function $f$ on $\s(A)\cup\s(B)$, where $C$ is an absolute constant.
\end{thm}

\medskip

{\bf Remark.}
The estimate in Theorem \ref{Cnr} is also sharp. Indeed, for each positive integer $n$ there exist $n\times n$ self-adjoint matrices $A$ and $B$ such that $A\ne B$ and
$$
\big\|\,|A|-|B|\,\big\|\ge \const\log(1+n)\|A-B\|,
$$
This follows easily from \rf{MG}, see the proof of Theorem 10.1 in \cite{AP2}.

\medskip

{\bf Definition.}
Let $\fF$ be a nonempty compact subset of $\R$. Recall that for $\e>0$, the $\e$-{\it entropy
$K_\e(\fF)$ of} $\fF$ is defined as
$$
K_\e(\fF)\df\inf\log\big(\card(\L)\big),
$$
where the infimum is taken over all $\L\subset\R$ such that $\L$ is an $\e$ net of $\fF$.
The following result is an generalization of Theorem \ref{skon}. On the other hand, it improves inequality \rf{modnep} obtained in \cite{AP2}.

\begin{thm}
\label{qcom}
Let $A$ and $B$ be self-adjoint operators and let $R$ be bounded operator
with $\|R\|\le1$. Suppose that $\s(A)\subset\fF$, where $\fF$ is a closed subset of $\R$. Then for every $f\in{\rm Lip}\big(\s(A)\cup\s(B)\big)$,
$$
\|f(A)R-Rf(B)\|\le\const\big(1+K_\e(\fF)\big)\|f\|_{{\rm Lip}(\s(A)\cup\s(B))}\|AR-RB\|,
$$
where $\e\df\|AR-RB\|$.
\end{thm}

\Pf We repeat the argument of the proof of Theorem \ref{kmr+}.
We can find a self-adjoint operator $A_\e$ such that
$A_\e A=AA_\e$, $\|A-A_\e\|\le\e$, $\s(A_\e)\subset\fF$,
and $\log\big(\card\big(\s(A_\e)\big)\big)\le K_\e(\fF)$. Then
\begin{align*}
\|f(A_\e)R-Rf(B)\|&\le\const\big(1+K_\e(\fF)\big)\|f\|_{{\rm Lip}(\s(A)\cup\s(B))}\|A_\e R-RB\|\\[.2cm]
&\le2\const\d\big(1+K_\e(\fF)\big)\|f\|_{{\rm Lip}(\s(A)\cup\s(B))}
\end{align*}
by Theorem \ref{skon}. It remains to observe that since $A$ commutes wit $A_\e$, we have
\begin{align*}
\|f(A)R-Rf(B)\|&\le\|f(A)-f(A_\e)\|+\|f(A_\e)R-Rf(B)\|\\[.2cm]
&\le\e\|f\|_{{\rm Lip}(\s(A))}+\|f(A_\e)R-Rf(B)\|.\quad \bl
\end{align*}

\begin{cor}
Let $A$ and $B$ be self-adjoint operators and let $\s(A)\subset\fF$, where $\fF$ is a closed subset of $\R$. Then
for every $f\in{\rm Lip}\big(\s(A)\cup\s(B)\big)$,
$$
\|f(A)-f(B)\|\le\const\big(1+K_\e(\fF)\big)\|f\|_{{\rm Lip}(\s(A)\cup\s(B))}\|A-B\|,
$$
where $\e\df\|A-B\|$.
\end{cor}

\Pf It suffices to put $R=I$. $\bl$

If we apply Theorem \ref{qcom} to the case $K=[a,b]$, we obtain the following estimate, which improves inequality \rf{modnep} in the special case $R=I$.

\begin{cor}
Let $f\in{\rm Lip}(\R)$. Let $A$ be a self-adjoint operator with $\s(A)\subset [a,b]$.
and $\|R\|\le1$. Then for every self-adjoint operators $B$,
$$
\|f(A)R-Rf(B)\|\le\const\|f\|_{\Li}\log\left(2+\frac{b-a}{\|AR-RB\|}\right)\|AR-RB\|.
$$
\end{cor}

\begin{cor}
Let $f\in{\rm Lip}(\R)$. Let $A$ be a self-adjoint operator with $\s(A)\subset [a,b]$.
and $\|R\|\le1$. Then
$$
\|f(A)R-Rf(B)\|\le\const\|f\|_{\Li}\log\left(2+\frac{b-a}{\|A-B\|}\right)\|A-B\|.
$$
\end{cor}

\

\

\noindent
\begin{tabular}{p{9cm}p{15cm}}
A.B. Aleksandrov & V.V. Peller \\
St-Petersburg Branch & Department of Mathematics \\
Steklov Institute of Mathematics  & Michigan State University \\
Fontanka 27, 191023 St-Petersburg & East Lansing, Michigan 48824\\
Russia&USA
\end{tabular}

\end{document}

%% file: classstartscr.tex
\newcommand{\lb}{\linebreak}

\renewcommand{\a}{\alpha}
\renewcommand{\b}{\beta}

\renewcommand{\d}{\delta}
\newcommand{\e}{\varepsilon}
\newcommand{\vk}{\varkappa}
\newcommand{\z}{\zeta}

\renewcommand{\l}{\lambda}

\newcommand{\s}{\sigma}
\renewcommand{\t}{\tau}
\newcommand{\var}{\varphi}
\newcommand{\f}{\varphi}
\renewcommand{\o}{\omega}

\newcommand{\D}{\Delta}
\renewcommand{\L}{\Lambda}

\renewcommand{\O}{\Omega}

\newcommand{\cc}{{\mathscr C}}
\newcommand{\E}{{\mathscr E}}

\newcommand{\F}{{\mathscr F}}

\newcommand{\cL}{{\mathscr L}}
\newcommand{\M}{{\mathscr M}}

\newcommand{\cp}{{\mathscr P}}

\newcommand{\C}{{\Bbb C}}
\newcommand{\T}{{\Bbb T}}

\newcommand{\dd}{{\Bbb D}}
\newcommand{\R}{{\Bbb R}}
\newcommand{\Z}{{\Bbb Z}}

\newcommand{\0}{{\boldsymbol{0}}}

\newcommand{\bs}{\boldsymbol}

\newcommand{\m}{{\boldsymbol m}}
\newcommand{\bS}{{\boldsymbol S}}

\newcommand{\rf}[1]{(\ref{#1})}

\newcommand{\df}{\stackrel{\mathrm{def}}{=}}
\newcommand{\dist}{\operatorname{dist}}

\newcommand{\supp}{\operatorname{supp}}

\newcommand{\rank}{\operatorname{rank}}
\newcommand{\const}{\operatorname{const}}

\newcommand{\eeq}{\end{equation}}
\newcommand{\beq}{\begin{equation}}
\newcommand{\bay}{\begin{eqnarray}}
\newcommand{\ba}{\begin{align*}}
\newcommand{\ea}{\end{align*}}
\newcommand{\ey}{\end{eqnarray}}
\newcommand{\bey}{\begin{eqnarray*}}
\newcommand{\eey}{\end{eqnarray*}}

\newcommand{\be}{\infty}

\newcommand{\bl}{\blacksquare}

\newcommand{\Pf}{{\bf Proof. }}

\newcommand{\ov}{\overline}

\newtheorem{thm}{\hspace{\parindent}Theorem}[section]

\newtheorem{cor}[thm]{\hspace{\parindent}Corollary}
\newtheorem{lem}[thm]{\hspace{\parindent}Lemma}